\documentclass[12pt, oneside,reqno]{amsart}

\usepackage{latexsym,amsxtra,amscd,ifthen}
\usepackage{amsfonts}
\usepackage{verbatim}
\usepackage{dsfont}
\usepackage{mathtools}
\usepackage{enumerate}
\usepackage[capitalise]{cleveref}
\usepackage{endnotes}

\usepackage{graphicx}
\usepackage[latin1]{inputenc}
\usepackage{array}
\usepackage{amsmath}
\usepackage{fancyhdr}
\usepackage{graphics}
\usepackage{amssymb}
\usepackage{amsthm}
\usepackage{enumitem}
\usepackage{mathrsfs}
\usepackage{yfonts} 
\usepackage{cancel}
\usepackage[square,numbers]{natbib}
\usepackage{stmaryrd}
\SetSymbolFont{stmry}{bold}{U}{stmry}{m}{n}
\usepackage{amsfonts} 
\usepackage[pdftex]{color}
\usepackage{ulem}
\usepackage{tikz-cd}		
\usepackage[arrow, matrix, curve]{xy}

\makeindex

\begin{document}

\def\mathllap{\mathpalette\mathllapinternal}
\def\mathllapinternal#1#2{%
\llap{$\mathsurround=0pt#1{#2}$}
}
\def\clap#1{\hbox to 0pt{\hss#1\hss}}
\def\mathclap{\mathpalette\mathclapinternal}
\def\mathclapinternal#1#2{%
\clap{$\mathsurround=0pt#1{#2}$}%
}
\def\mathrlap{\mathpalette\mathrlapinternal}
\def\mathrlapinternal#1#2{%
\rlap{$\mathsurround=0pt#1{#2}$}
}

\def\itemMath#1{\raisebox{-\abovedisplayshortskip}
{\parbox{1.0\linewidth}{\begin{equation}\begin{split}#1\end{split}\end{equation}}}}

\def\itemMathtwo#1#2{\raisebox{-#1pt}
{\parbox{1.0\linewidth}{\begin{equation}\begin{split}#2\end{split}\end{equation}}}}

\newenvironment{other}[1]{\noindent\textbf{#1 \hspace{0.1em} \stepcounter{thm}\thethm}\hspace*{1,5em}}{}

\let\amsamp=&

\newcommand{\CC}{\mathbb{C}}
\newcommand{\NN}{\mathbb{N}}
\newcommand{\QQ}{\mathbb{Q}}
\newcommand{\ZZ}{\mathbb{Z}}
\newcommand{\BBB}{\mathcal{B}}

\newcommand{\RR}{\mathbb{R}}
\newcommand{\PP}{\mathbb{P}}
\newcommand{\R}{\Rightarrow}
\newcommand{\p}{\varphi}
\newcommand{\nt}{\trianglelefteq}
\newcommand{\id}{\mathrm{id}}
\newcommand{\vol}{\operatorname{vol}}
\newcommand{\rank}{\operatorname{rank}}
\newcommand{\ann}{\operatorname{ann}}
\newcommand{\lcm}{\operatorname{lcm}}
\newcommand{\Sym}{\operatorname{Sym}}
\newcommand{\grad}{\operatorname{grad}}
\newcommand{\Aut}{\operatorname{Aut}}
\newcommand{\End}{\operatorname{End}}
\newcommand{\Hom}{\operatorname{Hom}}
\newcommand{\GrHom}{\operatorname{\underline{Hom}}}
\newcommand{\Mod}{\operatorname{Mod}}
\newcommand{\GrMod}{\operatorname{\underline{Mod}}}

\newcommand{\HS}{\operatorname{H}}
\newcommand{\Tr}{\operatorname{Tr}}
\newcommand{\BTr}{\operatorname{Br}}
\newcommand{\TrV}{\underline{\operatorname{Tr}}}
\newcommand{\TrVL}{\overline{\operatorname{Tr}}}
\newcommand{\HM}{\underline{\operatorname{H}}}
\newcommand{\tr}{\operatorname{tr}}
\newcommand{\gr}{\operatorname{gr}}
\newcommand{\op}{\operatorname{op}}
\newcommand{\hdet}{\operatorname{hdet}}
\newcommand{\Hdet}{\operatorname{Hdet}}
\newcommand{\grA}{\operatorname{gr}-\operatorname{A}}
\newcommand{\Ext}{\operatorname{Ext}}
\newcommand{\GrExt}{\operatorname{\underline{Ext}}}
\newcommand{\GK}{\operatorname{GKdim}}
\newcommand{\GKt}{\operatorname{GK}}
\newcommand{\gldim}{\operatorname{gl.dim}}
\newcommand{\pdim}{\operatorname{pdim}}
\newcommand{\injdim}{\operatorname{injdim}}
\newcommand{\proj}{\operatorname{proj}\,}
\newcommand{\dist}{\operatorname{dist}}
\newcommand{\cd}{\operatorname{cd}}
\newcommand{\qgr}{\operatorname{qgr}}
\newcommand{\tors}{\operatorname{tors}}
\newcommand{\sgn}{\operatorname{sgn}}
\newcommand{\type}{\operatorname{type}}
\newcommand{\Tor}{\operatorname{Tor}}
\newcommand{\adj}{\operatorname{adj}}
\newcommand{\col}{\operatorname{col}}

\newtheorem{innercustomgeneric}{\customgenericname}
\providecommand{\customgenericname}{}
\newcommand{\newcustomtheorem}[2]{%
  \newenvironment{#1}[1]
  {%
   \renewcommand\customgenericname{#2}%
   \renewcommand\theinnercustomgeneric{##1}%
   \innercustomgeneric
  }
  {\endinnercustomgeneric}
}


\numberwithin{equation}{section}
\newtheorem{theorem}[equation]{Theorem}
\newtheorem{lemma}[equation]{Lemma}
\newtheorem{proposition}[equation]{Proposition}
\newtheorem{corollary}[equation]{Corollary}
\newtheorem*{theorem*}{Theorem}

\theoremstyle{definition}
\newtheorem{convention}[equation]{Convention}
\newtheorem{defremark}[equation]{Definition \& Remark}
\newtheorem{definition}[equation]{Definition}
\newtheorem{example}[equation]{Example}

\theoremstyle{remark}
\newtheorem{remark}[equation]{Remark}
\newenvironment{thmbis}[1]
  {\renewcommand{\thetheorem}{\ref{#1}$'$}%
   \addtocounter{theorem}{-1}%
   \begin{theorem}}
  {\end{theorem}}








\title{Generalized Gorensteinness and a homological determinant for preprojective algebras}
\author{Stephan Weispfenning\protect\endnotemark[1]}

\subjclass[2010]{Primary: 16W22, 16W50, 16E65}

\keywords{invariant theory, preprojective algebras, homological determinant, generalized Gorenstein, trace function, Hilbert series, twisted Calabi-Yau}

\address{1. To whom correspondence should be addressed\hfill\break E-mail: sweispfe@ucsd.edu\hfill\break The author declares no conflict of interest.}
	
\begin{abstract}
The study of invariants of group actions on commutative polynomial rings has motivated many developments in commutative algebra and algebraic geometry. It has been of particular interest to understand what conditions on the group result in an invariant ring satisfying useful properties. In particular, Watanabe's Theorem states that the invariant subring of $k[x_1,\ldots,x_n]$ under the natural action of a finite subgroup of $SL_n(k)$ is always Gorenstein. In this paper, we study this question in the more general setting of group actions on noncommutative non-connected algebras $A$. We develop the notion of a homological determinant of an automorphism of $A$, then use the homological determinant to study actions of finite groups $G$ on $A$. We give a sufficient condition so that the invariant ring $A^G$ has finite injective dimension and satisfies the generalized Gorenstein condition. More precisely, let $A$ be a noetherian $\NN$-graded generalized Gorenstein algebra with finite global dimension. Suppose all elements of $G$ fix the idempotents of $A$ and act with trivial homological determinant. Then the invariant ring $A^G$ is generalized Gorenstein.
\end{abstract}

\maketitle

\pagestyle{fancy}
\fancyhf{}
\pagenumbering{arabic}\cfoot{\thepage}

\setlength{\headheight}{18pt}

\fancyhead[L]{\ifthenelse{\isodd{\value{page}}}{\center\small \tiny\rightmark}{\center \tiny STEPHAN WEISPFENNING}}

\section{Introduction}
In 1987 M.~Artin and W.~Schelter introduced a new class of noncommutative graded algebras which can be understood as noncommutative generalizations of polynomial rings: an algebra $S$ is \textit{Artin-Schelter regular} if it is a finitely generated connected $\NN$-graded ring $S = k \oplus \bigoplus_{j > 0}{S_j}$ satisfying the following three conditions:
\begin{itemize}[topsep=0pt, partopsep=0pt,leftmargin=0.75in]
	\item[(AS 1)] $S$ has finite (graded) global dimension $d \in \ZZ_{\geq 0}$. \par
	\item[(AS 2)] $S$ has finite Gelfand-Kirillov dimension. \par
	\item[(AS 3)] $S$ satisfies the Artin-Schelter Gorenstein condition
		\[
			\GrExt_S^m(k,S) \cong \GrExt_{S^{\op}}^m(k,S) \cong \delta_{md}k[\ell] \ \text{ \ \ for some } \ell \in \ZZ.
		\]
\end{itemize}
This class of noncommutative algebras has motivated a wide range of research interests. Many results in commutative invariant theory have been generalized to the noncommutative setting, with Artin-Schelter regular algebras playing the role of commutative polynomial rings. There has been extensive success in studying actions of groups (\cite{MR2434290, MR2601538, MR1758250}) and finite-dimensional semisimple Hopf algebras (\cite{MR2568355, MR3552496}) on Artin-Schelter regular algebras. 
For example, there are noncommutative analogues to the Chevalley-Shephard-Todd Theorem (\cite[Theorem 5.5]{MR2601538}), Watanabe's Theorem (\cite[Theorem 3.3]{MR1758250}, Stanley's Theorem (\cite[Proposition 3.8]{MR2568355}), and many more.\\ \indent
Of crucial importance for the results in this paper are the foundations laid by J{\o}rgensen and Zhang in \cite{MR1758250}.
More precisely, they prove the following remarkable theorem:
\begin{theorem*}[{\cite[Theorem 3.3]{MR1758250}}]
Suppose that $S$ is noetherian and AS-Gorenstein (i.e. $S$ has finite injective dimension over itself and satisfies (AS 3)), and that $G$ is a finite subgroup of $\Aut_{\gr}(S)$. If $\hdet(g) =1$ for each $g \in G$, then the fixed ring $S^G$ is AS-Gorenstein.
\end{theorem*}

In this paper we generalize this theorem to a broader class of so-called \textit{generalized Gorenstein algebras} which were first introduced in \cite[Definition 3.1]{MR2770441}. Let $B$ be a noetherian $\NN$-graded algebra with $B_0 \cong k^n$. Denote the primitive orthogonal idempotents of $B$ by $\{e_1,\ldots,e_n\}$, the simple right $B$-modules by $S_j = e_jB/(e_jB)_{\geq 1}$ and analogously the simple left $B$-modules by $S_j^\vee = Be_j/(Be_j)_{\geq 1}$. Then $B$ is called generalized Gorenstein if it satisfies
\begin{itemize}[topsep=0pt, partopsep=0pt]
	\item[(1)]	$\injdim_B(B) = N < \infty$ and
	\item[(2)]	the generalized Gorenstein condition
	\[
		\GrExt_B^i(S_j,B) \cong \delta_{Ni}S_{\sigma(j)}^\vee[\ell_j]
	\]
for some $\sigma \in \Sym(n)$ and some $\ell_j \in \ZZ$ for $j = 1,\ldots,n$.
\end{itemize}

Another important class of noncommutative algebras are twisted Calabi-Yau algebras. We extract the following definition from \cite[Definition 1.2]{DanToDo}. Let $B$ be a $k$-algebra and let $B^e$ denote the enveloping algebra $B \otimes_k B^{\op}$. We say that $B$ is \textit{twisted Calabi-Yau of dimension $d \geq 0$} if
\begin{itemize}[topsep=0pt, partopsep=0pt,leftmargin=0.75in]
	\item[(CY 1)] $B$ has a resolution of finite length by finitely generated projective $(B,B)$-bimodules and
	\item[(CY 2)] there exists an invertible $k$-central $(B,B)$-bimodule $U$ such that
	\[
		\Ext_{B^e}^i(B,B^e) \cong \begin{cases} 0, & i \neq d \\ U, & i = d \end{cases}
	\]
	as $(B,B)$-bimodules, where each $\Ext_{B^e}^i(B,B^e)$ is considered as a right $B^e$-module via the ``inner" right $B^e$-structure of $B^e$.
\end{itemize}

Twisted Calabi-Yau algebras generalize Artin-Schelter regular algebras (see \cite[Theorem 5.15]{DanToDo}). Indeed, by \cite[Lemma 1.2]{MR3250287}, a connected graded algebra is twisted Calabi-Yau if and only if it is Artin-Schelter regular. Further, Reyes and Rogalski show that the twisted Calabi-Yau condition is equivalent to the generalized Gorenstein condition \cite[Theorems 5.2 and 5.15]{DanToDo}.\\ \indent
One class of twisted Calabi-Yau algebras are the preprojective algebras of extended Dynkin quivers of type $A$, $D$, and $E$ which in addition have a natural grading, global dimension and Gelfand-Kirillov dimension equal to $2$ and nice homological properties. To define these algebras, let $Q$ be a quiver with finitely many vertices $Q_0 = \{i \mid 1 \leq i \leq n\}$ and finitely many arrows $Q_1$. The double $\bar{Q}$ of $Q$ is the quiver obtained by keeping the vertex set and adding a new arrow $\alpha^\ast$ from $j$ to $i$ for each arrow $\alpha$ from $i$ to $j$. The \textit{preprojective algebra} $\Pi(Q)$ is then defined as 
\[
	A = \Pi(Q) = k\bar{Q}/\left(\sum_{\alpha \in Q_1}{\alpha \alpha^\ast - \alpha^\ast \alpha}\right)
\]
where $k\bar{Q}$ denotes the path algebra of $\bar{Q}$. \\ \indent
In Section \ref{section6} we introduce a matrix homological determinant for a generalized Gorenstein algebra $B$ with $B_0 \cong k^n$ and $\injdim_B(B) = N$. This construction is a generalization of the homological determinant of Artin-Schelter regular algebras used in \cite[Section 2]{MR1758250}. We consider a graded automorphism $g$ of $B$ which fixes the primitive orthogonal idempotents $\{e_1,\ldots,e_n\}$. Lemma \ref{lemma65} connects $H_{\mathfrak{r}}^N(g)$ to the map $(g^{-1})^\ast$ where
\[
	H^i_{\mathfrak{r}}(M) = \varinjlim{\GrExt^i_B({}_B(B/B_{\geq r})_B,M_B)}
\]
for all right $B$-modules $M$, $\mathfrak{r} = B_{\geq 1}$ and 
\[
  (g^{-1})^\ast: \Hom_k(_BB_k,k_k) \to \Hom_k(_BB_k,k_k), \ \varphi \mapsto \varphi \circ g^{-1}. 
\]
More precisely, $H_{\mathfrak{r}}^N(g)$ equals $(g^{-1})^\ast$ up to $n$ nonzero scalars, one for each idempotent $e_i$. We write $\Hdet(g)$ for the matrix whose diagonal entries are the inverses of these scalars, and call it the homological determinant of $g$. This allows us to prove the following theorem.
\setcounter{section}{3}
\setcounter{equation}{11}	
\begin{theorem}
Let $B$ be a noetherian $\NN$-graded generalized Gorenstein algebra with degree zero piece $B_0 \cong k^n$. Suppose $B$ has finite global dimension $d$. Let $G$ be a finite subgroup of graded automorphisms of $B$ such that every $g \in G$ fixes the primitive idempotents. Assume that every $g \in G$ satisfies $\Hdet(g) = I_n$. Then $B^G$ is generalized Gorenstein.
\end{theorem}

In the final section, Section \ref{section7}, we focus on preprojective algebras $A$ of extended Dynkin quivers of type $A$, $D$ and $E$. It is well-known that in the connected case, the only Artin-Schelter regular algebras of global dimension $2$ are of the form $k \langle x, y \rangle /(r)$ for $r = yx - qxy$ for some nonzero $q \in k$ or $r = yx - xy - x^2$. In these cases, $g$ acts via scalar multiplication by $\hdet(g)$ on the defining relation $r$ (see \cite[Theorem 2.1]{MR3552496}). In the same way, we are able to conclude this paper with the following corollary.
\setcounter{section}{5}
\setcounter{equation}{7}	
\begin{corollary}
Let $Q$ be the quiver corresponding to an extended Dynkin graph of type $\widetilde{A_{n-1}}$ for $n \geq 2$, $\widetilde{D_{n-1}}$ for $n \geq 5$ or $\widetilde{E_m}$ for $m = 6,7,8$ and $A = \Pi(Q)$. Let $g$ be a graded automorphism which scales the arrows. Then $\Hdet(g) = c_1t_1 I_n$ ($n = m+1$ for $\widetilde{E_m}$).
\end{corollary}

\section*{Acknowledgments}
The author is extremely grateful for the support, help and proofreading of Daniel Rogalski. The results would not exist without his valuable insight and helpful discussions. The author would like to thank Rob Won for proofreading parts of this paper. The author was partially supported by NSF grants DMS-1201572 and DMS-1601920 as well as NSA grant H98230-15-1-0317. The author also wishes to thank Efim Zelmanov for the generous financial support that allowed him to fully focus on this project.

\setcounter{section}{1}
\setcounter{equation}{0}	

\section{Preliminaries}
Throughout let $k$ be an algebraically closed field of characteristic $0$. Unless stated otherwise, every algebra is $\NN$-graded and locally finite, i.e., each graded piece is finite dimensional over $k$. All quivers have finitely many vertices and finitely many arrows. The symbol $A$ is reserved to be the preprojective algebra $\Pi(Q)$ of an extended Dynkin quiver $Q$ (see Euclidean graphs in \cite[Chapter VII.2]{MR2197389}) as defined in Definition \ref{def23}. \\ \indent
We start by fixing our notation. Let $R = \bigoplus\limits_{j \in \NN}{R_j}$ be a graded algebra over $k$ and let $g$ be an automorphism of $R$.

\begin{itemize}[topsep=0pt]
	\item	We call $g$ a \textit{graded automorphism} if $g(R_i) = R_i$ for all graded pieces $R_i$ of $R$. The group of all graded automorphisms of $R$ is denoted $\Aut_{\gr}(R)$.
	\item Let $G$ be a subgroup of $\Aut_{\gr}(R)$. The \textit{invariant ring} or \textit{fixed ring} is defined as $R^G = \{x \in R \mid g(x) = x, \text{ for all } g \in G\}$.
	\item	A right $R$-module $M$ is called \textit{graded} if $M = \bigoplus_{j \in \ZZ}{M_j}$ and all $r \in R_n$ and $m \in M_\ell$ satisfy $m r \in M_{n+\ell}$ for all $n, \ell \in \ZZ$. The category of graded right $R$-modules is denoted by $\GrMod$-$R$.
	\item	A graded right $R$-module $M$ is called \textit{left bounded} if there exists an integer $N$ such that $M_j = 0$ for all $j \leq N$. Analogously, $M$ is called \textit{right bounded} if there exists an integer $N$ such that $M_j = 0$ for all $j \geq N$. The module $M$ is called \textit{bounded} if it is both left bounded and right bounded. Notice that every finitely generated $R$-module is left bounded.
	\item	For a graded $R$-module $M = \bigoplus_{j \in \ZZ}{M_j}$, the \textit{shifted module} $M[\ell] = \bigoplus_{j \in \ZZ}{(M[\ell])_j}$ is defined by $(M[\ell])_j = M_{\ell+j}$.
	\item	For graded $R$-modules $M = \bigoplus_{j \in \ZZ}{M_j}$ and $N = \bigoplus_{i \in \ZZ}{N_i}$, a linear map $\varphi$ satisfying $\varphi(M_j) \subseteq N_{j+\ell}$ is called \textit{homogeneous of degree $\ell$}. Denote the vector space of all right $R$-module homomorphisms from $M$ to $N$ which are homogeneous of degree $\ell$ by $\Hom^\ell_R(M,N)$.
	\item	For graded $R$-modules $M$ and $N$, the elements of $\Hom^0_R(M,N)$ are called \textit{graded} and the set of all graded automorphisms of $M$ is denoted by $\Aut_{\gr}(M)$.
	\item	For graded $R$-modules $M$ and $N$, the vector space $\Hom_R(M,N)$ denotes the usual $\Hom$ in the category $\Mod$-$R$ of all right $R$-modules. We have
	\[
		\GrHom_R(M,N) = \bigoplus_{\ell = -\infty}^{\infty}{\Hom_R^\ell(M,N)} \subseteq \Hom_R(M,N).
	\]
	It is well-known that we have equality whenever $M$ is finitely generated which is the case for all appearances of $\GrHom$ in this paper. 
	\item	Denote the right derived functors of $\GrHom$ by $\GrExt_R^i = R^i\GrHom_R$.
	\item	For a right $R$-module $M$, the \textit{$g$-twisted module} $M^g$ is defined as $M^g = M$ with changed action $m \ast r = m g(r)$ for $m \in M, \ r \in R$.
	\item	Let $M$ and $N$ be graded right $R$-modules. A $k$-linear map $\varphi: M \to N$ is called \textit{$g$-linear} if $\varphi(m r) = \varphi(m) g(r)$ for all $m \in M$ and $r \in R$. A graded map $\varphi: M \to N$ is $g$-linear if and only if it is a graded $R$-homomorphism from $M$ to $N^g$.
\end{itemize}

\begin{definition}
A \textit{quiver} $Q = (Q_0,Q_1,s,t)$ is a quadruple consisting of vertices $Q_0 = I$ for some index set $I$, arrows $Q_1$ and two maps $s,t: Q_1 \to Q_0$ called the \textit{source} $s$ and \textit{target} $t$. For every arrow $\alpha \in Q_1$ the maps associate the corresponding source $s(\alpha)$ and target vertex $t(\alpha)$ in $Q_0$. Usually, the quiver is just denoted by $Q$. 
A \textit{path of length} $\ell \geq 1$ with source $a$ and target $b$ is a sequence
\[
	\alpha_1 \alpha_2 \cdots \alpha_\ell 
\]
where $\alpha_j \in Q_1$ for all $1 \leq j \leq \ell$ as well as $s(\alpha_1) = a, \ t(\alpha_j) = s(\alpha_{j+1})$ for each $1 \leq j < \ell$ and $t(\alpha_\ell) = b$. Additionally, each vertex $i \in Q_0$ has a path of length $\ell = 0$ associated to it, called the \textit{trivial path} at $i$ denoted by $e_i$.\\ \indent
The \textit{path algebra} $kQ$ of $Q$ is the $k$-algebra whose underlying vector space has as its basis the set of all paths of length $\ell \geq 0$ such that the product of two paths $\alpha_1 \alpha_2 \cdots \alpha_\ell$ and $\beta_1 \beta_2 \cdots \beta_m$ is defined by
\[
\alpha_1 \alpha_2 \cdots \alpha_\ell \cdot \beta_1 \beta_2 \cdots \beta_m = \delta_{t(\alpha_\ell)s(\beta_1)} \alpha_1 \cdots \alpha_\ell \beta_1 \cdots \beta_m
\]
where $\delta_{t(\alpha_\ell)s(\beta_1)}$ denotes the Kronecker delta. The path algebra has a natural grading 
\[
	kQ = kQ_0 \oplus kQ_1 \oplus kQ_2 \oplus \ldots \oplus kQ_{\ell} \oplus \ldots
\]
where $kQ_\ell$ denotes the $k$-subspace of $kQ$ generated by the set $Q_\ell$ of all paths of length $\ell$.
\end{definition}

\begin{definition}
For a quiver $Q$, its \textit{double} $\bar{Q}$ is obtained by keeping the vertex set $Q_0$ and adding a new arrow $\alpha^\ast$ from $j$ to $i$ for each arrow $\alpha \in Q_1$ from $i$ to $j$.
\end{definition}

\begin{definition}\label{def23}
Using \cite[Definition 3.1.1.]{MR2335985}, we can define the preprojective algebra associated to a quiver $Q$. Let $Q$ be a quiver and $\bar{Q}$ its double. The \textit{preprojective algebra} is defined as
\[
	A = \Pi(Q) = k\bar{Q}/\left(\sum_{\alpha \in Q_1}{\alpha \alpha^\ast - \alpha^\ast \alpha}\right).
\]
\end{definition}

Note that we are using the version of the definition in which loops are doubled. A weak homological property introduced by Artin and Zhang is the $\chi$-condition for a noetherian $\NN$-graded algebra from \cite[Definitions 3.2]{MR1304753}.  
\begin{definition}\label{defchi}
A noetherian $\NN$-graded algebra $R$ satisfies the \textit{$\chi$-condition} if $\GrExt_R^j(R_0,M)$ is bounded for all $j$ and all finitely generated right $R$-modules $M$. 
\end{definition}

Our definition for $\chi$ is not the original definition given in \cite[Definition 3.7]{MR1304753}. However, due to \cite[Proposition 3.11 (2)]{MR1304753} the original definition is equivalent to Definition \ref{defchi} as long as $R$ is locally finite, which is assumed to be the case for all our algebras. We now recall some results about preprojective algebras of extended Dynkin quivers. 

\begin{proposition}\label{prop25}	Let $A = \Pi(Q)$ where $Q$ is an extended Dynkin quiver with vertex set $\{e_1,\ldots,e_n\}$. We denote the simple right $A$-modules by $S_j = e_jA/(e_jA)_{\geq 1}$ and analogously the simple left $A$-modules by $S_j^\vee = Ae_j/(Ae_j)_{\geq 1}$. 
 \begin{itemize}[topsep=0pt]
	\item[(a)]	{\normalfont\cite[Proposition (2.11)]{OtherPaper}}. $A$ has global dimension $2$ and satisfies the generalized Gorenstein condition 
	\[
		\GrExt^i_A(S_j,A) \cong \delta_{di} S_{\sigma(j)}^\vee[\ell_j]
	\]
	for all $j = 1,\ldots,n$, $\sigma = \id$ and $\ell_j = 2$ for all $j$.
	\item[(b)]	{\normalfont\cite[Proposition 5.1]{MR1304753}}. $A$ satisfies the $\chi$-condition from Definition \ref{defchi}.
	\item[(c)]	{\normalfont\cite[Theorem]{MR624903}}. $\GK(\Pi(Q)) = 2$. (The typo was noted in {\normalfont\cite[Proposition 6.10]{MR876985}}).
	\item[(d)]	{\normalfont\cite[Section 3]{DanToDo2}}. The trivial paths $\{e_1,\ldots,e_n\}$ associated to the vertices form the set of primitive and orthogonal idempotents. Then, the modules $P_i = e_iA$ for $i = 1,\ldots,n$ are the only indecomposable projectives and thus every finitely generated graded projective module is a sum of finitely many possibly shifted copies of the $P_i$'s.
	\item[(e)]	{\normalfont\cite[Theorem 6.5]{MR876985}}. $A$ is a noetherian prime polynomial identity ring.
\end{itemize}
\end{proposition}

Since $S_j$ is a right $A$-module, the dual space $\Hom_k(_kS_j,{}_kk)$ is a simple left $A$-module which is isomorphic to $Ae_j/(Ae_j)_{\geq 1}$ and that justifies the notation $S_j^\vee$ for the simple left $A$-modules. \\ \indent
Most of the main theory developed in Section \ref{section6} is worked out for the more general class of algebras which are generalized Gorenstein. The following definition was introduced in \cite[Definition 3.1]{MR2770441}.

\begin{definition}\label{genGS}
A noetherian $\NN$-graded algebra $B$ with degree zero piece $B_0 \cong k^n$ and primitive orthogonal idempotents $\{e_1,\ldots,e_n\}$ is called \textit{generalized Gorenstein} if it satisfies
\begin{itemize}[topsep=0pt]
	\item[(1)]	$\injdim_B(B) = N < \infty$ and
	\item[(2)]	the \textit{generalized Gorenstein condition} 
	\[
		\GrExt_B^i(S_j,B) \cong \delta_{Ni}S_{\sigma(j)}^\vee[\ell_j]
	\]
	for some $\sigma \in \Sym(n)$ where $S_i = e_iB/(e_iB)_{\geq 1}$ and $S_i^\vee = Be_i/(Be_i)_{\geq 1}$ are the simple right and left $B$-modules and $\ell_j \in \ZZ$.
\end{itemize}
\end{definition}

A weaker condition is the AS-Cohen-Macaulay condition (see \cite[Definition 0.1]{MR1758250}).
\begin{definition}\label{ASCM}
An algebra $R$ is called \textit{AS-Cohen-Macaulay} if there exists an integer $N$ such that for $i \neq N$ it follows that $H_{\mathfrak{u}}^i(R) = H_{\mathfrak{u}^{\op}}^i(R) = 0$ where $H_{\mathfrak{u}}^i(R) = \varinjlim{\GrExt^i({}_R(R/\mathfrak{u}^n)_R,R_R)}$ and $\mathfrak{u} = R_{\geq 1}$.
\end{definition}

Let's recall the definitions of the Hilbert series and the trace function which play an important role in any kind of invariant theory.
\begin{definition}\label{def28}
Let $R$ be an $\NN$-graded algebra over $k$ and let $M$ be a left bounded graded locally finite right $R$-module. Further, let $g \in \Aut_{\gr}(M)$ be a graded automorphism of $M$. Then:
\begin{itemize}[topsep=0pt]
	\item	The \textit{Hilbert series} $\HS_M(t)$ of $M$ is the formal Laurent series
\[
	\HS_M(t) = \sum_{j \in \ZZ}{\dim_k(M_j)t^j}.
\]	
	\item	The \textit{trace} $\Tr_M(g,t)$ of $g$ is the formal Laurent series
\[
	\Tr_M(g,t) = \sum_{j \in \ZZ}{\tr_k(g|_{M_j})t^j}. 
\]
\end{itemize}
\end{definition}

\begin{definition}\label{def29}
Let $B$ be a noetherian $\NN$-graded algebra with $B_0 \cong k^n$. Denote the primitive orthogonal idempotents by $\{e_1,\ldots,e_n\}$. Let $g \in \Aut_{\gr}(B)$ be a graded automorphism of $B$ such that $g(e_i) = e_i$ for every $i = 1,\ldots,n$. Notice that $g$ restricts to an automorphism of $e_iB$, and that there is an induced map $g_i: e_iB \longrightarrow e_iB, \ e_i \beta \mapsto e_i g(\beta)$ which is $g$-linear.
\end{definition}

For non-connected algebras it is important to find a generalization of the trace function which encodes more information and allows us to consider the indecomposable projective modules one at a time. Therefore, we will introduce the concept of vector traces of a graded automorphism which fixes the idempotents.

\begin{definition}
Let $B$ be a noetherian $\NN$-graded algebra with $B_0 \cong k^n$. Denote the primitive orthogonal idempotents by $\{e_1,\ldots,e_n\}$. For $g \in \Aut_{\gr}(B)$ with $g(e_i) = e_i$, the \textit{vector traces} $\TrV_B(g,t)$ and $\TrVL_B(g,t)$ (or $\TrV(g,t)$ and $\TrVL(g,t)$, respectively) are the vectors with entries
\begin{align*}
	[\TrV_B(g,t)]_{j} &= [\TrV(g,t)]_{j} = \sum_{s=0}^{\infty}{\tr(g|_{e_jB_s})t^s}, \ \ \ \text{ and} \\
	[\TrVL_B(g,t)]_{j} &= [\TrVL(g,t)]_{j} = \sum_{s=0}^{\infty}{\tr(g|_{B_se_j})t^s}.
\end{align*}
\end{definition}

The remaining results of this section work in very general settings. Lemma \ref{lemma211} only requires a ring $B$ to be graded and provides a very interesting projection of $B$ onto the fixed ring under a finite group action. More importantly, it makes $B$ into a finitely generated module over its fixed ring. This is one of the key ingredients to connect $H_{\mathfrak{r}}^N(B)$ to $H_{\mathfrak{s}}^N(B^G)$ for $\mathfrak{r} = B_{\geq 1}$ and $\mathfrak{s} = (B^G)_{\geq 1}$ in Section \ref{section6}. 

\begin{lemma}\label{lemma211}	Let $R$ be an $\NN$-graded ring and $G$ be a finite subgroup of $\Aut_{\gr}(R)$. Then:
\begin{itemize}[topsep=0pt]
	\item[(a)]	{\normalfont \cite[Lemma 5.1]{MR1438180}}. $R \cong R^G \oplus K$ as $(R^G,R^G)$-bimodules, where $K$ is the kernel of the map
	\[
		\pi_G: x \mapsto \frac{1}{|G|}\sum_{g\in G}{g(x)}.
	\]
	\item[(b)]	{\normalfont \cite[Corollary 1.12]{MR590245}}. If $R$ is right noetherian, then $R^G$ is right noetherian.
	\vspace{18pt}
	\item[(c)]	Molien's Theorem, {\normalfont\cite[Lemma 5.2]{MR1438180}}. The Hilbert series of the fixed ring equals: 
	\[
		\HS_{R^G}(t) = \frac{1}{|G|} \sum_{g \in G}{\Tr_R(g,t)}.
	\]
\end{itemize}
\end{lemma}

We conclude this introduction by stating two important results. 
\begin{proposition}[{\cite[Corollary 5.9]{MR590245}}]\label{prop212}
Let $R$ be a right (left) noetherian ring, $G$ be a finite subgroup of $\Aut(R)$ such that $|G|^{-1} \in R$. Then $R$ is a finitely generated right (left) $R^G$-module.
\end{proposition}

\begin{proposition}[{\cite[Remark after Proposition 8.7]{MR1304753}}]\label{prop213}
Let $R$ be a noetherian graded ring and $G$ be a finite subgroup of $\Aut_{\gr}(R)$. If $R$ satisfies the $\chi$-condition then so does $R^G$.
\end{proposition}

\begin{convention}\label{con214}
Throughout let $k$ be an algebraically closed field of characteristic $0$. Every algebra is noetherian, $\NN$-graded with degree zero piece isomorphic to $k^n$ and locally finite. Our convention is to use right modules. Every automorphism obeys the grading and all groups $G$ and automorphisms $g$ have finite order. Moreover, every $g \in G$ fixes the idempotents. Unless said otherwise, $A = \Pi(Q)$ for an extended Dynkin quiver $Q$ of type $\widetilde{A_{n-1}}$ for $n \geq 2$, $\widetilde{D_{n-1}}$ for $n \geq 5$ or $\widetilde{E_m}$ for $m = 6,7,8$ with primitive orthogonal idempotents $\{e_1,\ldots,e_n\}$. Every automorphism $g$ of $A$ that we consider scales the arrows. In other words, for each pair of arrows $\alpha$ and $\alpha^\ast$, there exist nonzero scalars $c_\alpha$ and $t_\alpha$ such that $g(\alpha) = c_\alpha \alpha$ and $g(\alpha^\ast) = t_\alpha \alpha^\ast$ which fully describes the action of $g$. This restriction is equivalent to the standing condition that $g$ fixes the primitive orthogonal idempotents except for $Q = \widetilde{A_1}$ where it is slightly stronger. Moreover, the relations of $A$ force $c_\alpha t_\alpha = c_\beta t_\beta$ for all arrows $\alpha, \beta \in Q_1$. Unless said otherwise, $B$ denotes a noetherian $\NN$-graded generalized Gorenstein algebra with $B_0 \cong k^n$ such that $B$ and $B^{\op}$ satisfy the $\chi$-condition. To simplify notation, the primitive orthogonal idempotents of $B$ are also denoted by $\{e_1,\ldots,e_n\}$. The modules $S_i$ and $S_i^\vee$ denote the simple right and left modules of the algebra we focus on. We call $\mathfrak{m} = A_{\geq 1}, \ \mathfrak{r} = B_{\geq 1}$ and $\mathfrak{s} = (B^G)_{\geq 1}$. The invariant ring or fixed ring is denoted by $B^G = \{b \in B \mid g(b) = b, \text{ for all } g \in G\}$.
\end{convention}

\section{A sufficient condition for an invariant ring to be generalized Gorenstein}\label{section6}
Fix an algebra $B$, a finite subgroup $G$ of $\Aut_{\gr}(B)$ and some $g \in G$ as in Convention \ref{con214}.\\ \indent
In this section, we develop the main theory of this paper. In particular, we introduce the concept of a homological determinant $\Hdet(g)$ in matrix form for the graded automorphism $g$. The main theorem, Theorem \ref{thm69} proves that whenever the homological determinant $\Hdet(g)$ for all $g$ in $G$ is trivial then the fixed ring is generalized Gorenstein. While stated in a more general setting, the results of this section apply to preprojective algebras. \\ \indent
The following is a generalization of the work of Peter J{\o}rgensen and James J. Zhang in \cite{MR1758250}. Denote $B_{\geq 1}$ by $\mathfrak{r}$ and recall that for a right $B$-module $M$
\begin{align*}
	H^i_{\mathfrak{r}}(M) &= \varinjlim{\GrExt^i({}_B(B/B_{\geq r})_B,M_B)}
\end{align*}
has a right $B$-module structure since $B/B_{\geq r}$ is a $(B,B)$-bimodule. The next definition is adapted from \cite{MR1304753}:
\begin{definition}
An element $e$ of a $B$-module $M$ is called \textit{torsion} if there exists an $i \in \NN$ such that $eB_{\geq i} = 0$. 
\end{definition}

In order to avoid confusion with the arrows $\alpha^\ast$ of $A$, we denote elements of a dual space by $(-)^\vee$. In the same way, the dual space $\GrHom_k(B,k)$ of $B$ is denoted by $B^\vee$ which has a natural $(B,B)$-bimodule structure induced by the $(B,B)$-bimodule structure of $B$. Similarly, $\GrHom_k(Be_j,k)$ is naturally a right $B$-module induced by the left $B$-module structure of $Be_j$. By \cite[Lemma 2.8]{DanToDo} this encourages us to call 
\[
	\GrHom_k(Be_j,k) = e_j\GrHom_k(B,k) = e_jB^\vee
\]
as right $B$-modules. The following key lemma shows that this module equals the graded injective hull of the simple module $S_j$.

\begin{lemma}\label{lemma62}
The graded Matlis dual $e_jB^\vee = \GrHom_k(Be_j,k)$ of the left $B$-module $Be_j$ with natural grading $(e_jB^\vee)_m = \GrHom_k(B_{-m}e_j,k)$ is the graded injective envelope of $S_j$.
\begin{proof}
Fix some homogeneous $k$-basis $\BBB$ of $Be_j$ which contains $e_j$. First, the graded $B$-module homomorphism $\varphi: S_j \to e_jB^\vee, \ e_j \mapsto e_j^\vee$ where $e_j^\vee$ denotes the dual map to $e_j$, sending $e_j$ to $1$ and all other elements of $\BBB$ to $0$, gives an inclusion of $S_j$ into $e_jB^\vee$.\\ \indent
Second, we need to show that $e_jB^\vee$ is the graded injective hull. Equivalently, it is enough to prove that $e_jB^\vee$ is an injective essential extension of $S_j$. Let $N$ be a nonzero submodule of $e_jB^\vee$. We denote the dual basis of $e_jB^\vee$ by $\BBB^\vee = \{\beta^\vee \mid \beta \in \BBB\}$. Let $\psi = \sum_{\beta \in \BBB}{c_\beta \beta^\vee}: Be_j \to k$ be a nonzero element from $N$. Since $e_jB^\vee$ lives entirely in nonpositive degrees, let $\gamma^\vee$ be of minimal degree such that $c_\gamma$ is nonzero. We are going to show that $\psi \cdot \gamma \in N \cap S_j \setminus \{0\}$. For $\beta^\vee$ with $\deg(\beta^\vee) > \deg(\gamma^\vee)$ in the expression for $\psi$ we notice that
\[
	(\beta^\vee \cdot \gamma)(\alpha) = \beta^\vee(\gamma \alpha) = 0
\]
since $\gamma \alpha$ can never equal $\beta$ as its degree is too big. Also, elements $\beta^\vee \in \BBB^\vee$ with $\deg(\beta^\vee) = \deg(\gamma^\vee)$ satisfy $\beta^\vee \cdot \gamma = e_j^\vee$ if and only if $\beta = \gamma$ and $0$ otherwise. Combined, we obtain
\[
	\psi \cdot \gamma = c_\gamma e_j^\vee
\]
which lives in $N \cap S_j$. Hence, $e_jB^\vee$ is an essential extension.\\ \indent
In regards to injectivity, recall that $P_j$ is a projective locally finite right $B$-module. Let $M = \bigoplus_i M_i$ be a locally finite $B$-module. Then $\GrHom_k(\bigoplus M_i,k) = \bigoplus \GrHom_k(M_{-i},k)$ together with 
\[
	\GrHom_k(\GrHom_k(M_i,k),k) \cong M_i
\]
implies that taking the Matlis dual twice returns a module isomorphic to $M$. Therefore, the functor $\GrHom_k(-,k)$ forms a duality on the category of locally finite graded $B$-modules. In particular, the graded projective left $B$-module $Be_j$ is sent to a graded injective right $B$-module $e_jB^\vee$.
\end{proof}
\end{lemma}

This will now help proving the following lemma motivated by \cite[Lemma 2.1]{MR1758250}. Since $A$ is generalized Gorenstein by Proposition \ref{prop25} this lemma also applies to $A$.

\begin{lemma}\label{lemma63}
Recall that $B$ satisfies $\GrExt_B^i(S_j,B) \cong \delta_{Ni}S_{\sigma(j)}^\vee[\ell_j]$ for some $\sigma \in \Sym(n)$. For $\mathfrak{r} = B_{\geq 1}$ we have
\[
	H^i_{\mathfrak{r}}(B) \cong \begin{cases} 0 & \text{ for } i \neq N, \\ \bigoplus_{b=1}^{n}{e_bB^\vee}[\ell_b]  & \text{ for } i = N \end{cases}
\]
where $e_bB^\vee = \GrHom_k(_B(Be_b) _k,k_k)$ and $N = \injdim(B)$. In particular, $B$ is AS-Cohen-Macaulay.
\begin{proof}
We calculate $\GrExt^i(B/B_{\geq r},B)$ by resolving the second entry. For this, we use a minimal graded injective resolution of length $N$ of $B$, say
\[
\begin{tikzcd}
0 \arrow[r] & B \arrow[r] & E^0 \arrow[r] & E^1 \arrow[r] & \ldots  \arrow[r] & E^N \arrow[r] & 0.
\end{tikzcd}
\]
Let $j \leq N$ be minimal such that $E^j$ has a nonzero torsion element $e$. We can assume that $eB_{\geq i} = 0$ for some minimal $i \in \NN$. By picking a nonzero element $\tilde{e} \in eB_{\geq i-1}$ we have that $\tilde{e}B_{\geq 1} = 0$. Since $\tilde{e} \cdot e_a \neq 0$ for some $a$, there is an embedding $S_{a}[-\deg(\tilde{e})] \hookrightarrow E^j$. As $S_a$ is a simple module and the resolution is minimal, this forces the composition $S_a \to E^j \to E^{j+1}$ to be zero as otherwise the graded injective hull of $S_a$ is a summand of both $E^j$ and $E^{j+1}$. Thus, applying $\GrHom(S_a[-\deg(\tilde{e})],-)$ to the resolution and taking homology at the $j$-th spot, we obtain nonzero elements. In other words, using the minimality of $j$, $\GrExt_B^j(S_a[-\deg(\tilde{e})],B) \neq 0$ and so $j = N$ as well as $E^0,\ldots,E^{N-1}$ are torsionfree. This shows that $\GrHom(B/B_{\geq r},E^i) = 0$ for $i < N$ and we can conclude $H^i_{\mathfrak{r}}(B) = 0$ for $i \neq N$.\\ \indent
From the generalized Gorenstein condition it follows that $\GrExt_B^N(S_b,B) \cong S_{\sigma(b)}^\vee[\ell_b]$ for every $b$ and some $\sigma \in \Sym(n)$. As shown above, $E^0,\ldots,E^{N-1}$ are torsionfree. Fixing $b \in \{1,\ldots,n\}$, this means $\GrHom_B(S_b,E^N)$ is nonzero. More precisely, due to $S_{b}$ being simple, we can find a right $B$-module embedding from $S_{b}$ into $E^N$. Since $E^N$ is an injective module, we obtain that $S_b[\ell_b]$'s graded injective envelope $e_bB^\vee[\ell_b]$ (see Lemma \ref{lemma62}) is a direct summand in $E^N$:
\[
	E^N \cong e_bB^\vee[\ell_b] \oplus M.
\]
To guarantee that $\dim_k(\GrExt_B^N(S_b,B)) = 1$, there cannot be another copy of $e_bB^\vee$ in $M$. Repeating this argument for $b = 1,\ldots,n$, we obtain
\[
	E^N \cong \bigoplus_{b=1}^{n}{e_bB^\vee}[\ell_b] \oplus I
\]
where $I$ is torsionfree. In order to finish the proof we use that $e_bB^\vee$ is torsion. Taking $H_{\mathfrak{r}}^N(-)$ provides
\begin{align*}
	H_{\mathfrak{r}}^N(B) &= \varinjlim{\GrExt^N(B/B_{\geq r},B)} = \varinjlim{\GrHom_B(B/B_{\geq r},E^N)} \\
	&= \varinjlim{\GrHom_B\left(B/B_{\geq r} \ ,\ \bigoplus_{b=1}^{n}{e_bB^\vee}[\ell_b] \oplus I\right)} = \bigoplus_{b=1}^{n}{e_bB^\vee}[\ell_b]
\end{align*}
where the last equality is true since the direct limit of $\GrHom_B(B/B_{\geq r},M)$ always gives back the torsion part of $M$. This implies the statement.
\end{proof}
\end{lemma}
Lemma \ref{lemma63} shows that $H_{\mathfrak{r}}^i(B) \cong \Hom_k(B,k) = B^\vee$ as ungraded modules. In general, the shifts of $e_bB^\vee[\ell_b]$ can be different for each $b$. Preprojective algebras on the other hand behave nicely, each $\ell_b = 2$ and therefore $H_{\mathfrak{m}}^2(A) \cong A^\vee[2]$ as graded modules. \\ \indent
Let $g \in \Aut_{\gr}(B)$ with $g(e_i) = e_i$ for all $i = 1,\ldots,n$. With the lemma just proved, our goal will be to connect $H_{\mathfrak{r}}^N(g): \ H_{\mathfrak{r}}^N(B) \to H_{\mathfrak{r}}^N(B)$ to $(g^{-1})^\ast$. Therefore, we need to understand $(g^{-1})^\ast$ first. It turns out that any $g$-linear map $\Phi: B^\vee_B \to B^\vee_B$ is closely related to $(g^{-1})^\ast$ via $n$ scalars determined by $\Phi(e_i^\vee)$. Notice that it is enough to look at $B^\vee$ rather than $\bigoplus\limits_{b = 1}^{n}{e_bB^\vee[\ell_b]}$ since $\Phi$ is $g$-linear and $g$ fixes the idempotents. This guarantees that $\Phi$ restricts to maps on each $e_bB^\vee[\ell_b]$ individually. For simplicity, we will therefore state the following remark and its subsequent lemma without worrying about shifts in the grading.
\begin{remark}\label{remark64}
We recall that $g \in \Aut_{\gr}(B)$ with $g(e_i) = e_i$ for all $i$. Then $g^{-1}$ becomes a $g^{-1}$-linear $B$-module homomorphism from $_BB$ to $_BB$. Furthermore, it is easy to see that
\[
	(g^{-1})^\ast: \ B^\vee_B \to \ B^\vee_B, \ \ \ \ \varphi \mapsto \varphi \circ g^{-1}
\]
is $g$-linear. Any $g$-linear $B$-module homomorphism from $B^\vee_B \to B^\vee_B$, say $\Phi: B^\vee_B \to B^\vee_B$, can be described as follows. We denote the dual basis by $\BBB^\vee = \{\beta^\vee \mid \beta \in \BBB\}$ where $\BBB$ stands for the union of fixed bases of $e_rBe_s$ for all $1 \leq r, s \leq n$. In general, we have
\[
	\Phi(e_i^\vee) = \sum_{j = 1}^{n}{\lambda_{ij} e_j^\vee}, \ \ \ \ \ \ \ \ \ \ \lambda_{ij} \in k.
\]
However, since 
\[
	\sum_{j = 1}^{n}{\lambda_{ij} e_j^\vee} = \Phi(e_i^\vee) = \Phi(e_i^\vee \cdot e_i) = \Phi(e_i^\vee) \cdot g(e_i)  = \sum_{j = 1}^{n}{\lambda_{ij} e_j^\vee} \cdot g(e_i) = \lambda_{ii} e_i^\vee
\]
we can define $\lambda_i = \lambda_{ii}$ and know that $\Phi(e_i^\vee) = \lambda_i e_i^\vee.$ We will show that these $\lambda_i$ determine $\Phi$ completely. Let $\beta$ be a basis element in $Be_s$. By definition it follows that $\beta^\vee \cdot \beta = e_s^\vee$. Now, with the same trick as for $e_i^\vee$, we can calculate
\[
	\lambda_s e_s^\vee = \Phi(e_s^\vee) = \Phi(\beta^\vee \cdot \beta) = \Phi(\beta^\vee) \cdot g(\beta).
\]
By checking on a basis, we conclude with $\Phi(\beta^\vee) = \lambda_s (\beta^\vee \circ g^{-1})$ for $\beta \in Be_s$.
\end{remark}
In order to apply the previous remark it remains to show that $H_{\mathfrak{r}}^N(g)$ is actually $g$-linear. Thanks to Lemma \ref{lemma63}, the domain and codomain equal $B^\vee_B$ as ungraded modules.
\begin{lemma}\label{lemma65}
The map $H_{\mathfrak{r}}^N(g): H_{\mathfrak{r}}^N(B) \to H_{\mathfrak{r}}^N(B)$, where $\mathfrak{r} = B_{\geq 1}$, is a $g$-linear $B$-module homomorphism from $B^\vee_B$ to $B^\vee_B$. Further, there exist scalars $\lambda_i$ for $i = 1,\ldots,n$ such that $H_{\mathfrak{r}}^N(g)(\beta^\vee) = \lambda_s (\beta^\vee \circ g^{-1})$ for $\beta \in Be_s$.
\begin{proof}
Thanks to Lemma \ref{lemma63}, this statement makes sense. We adapt the ideas of \cite[Lemma 2.2]{MR1758250} to our situation. The same proof, considering the $B$-linear map $g: B \to B^g$ (the right $B$-module with action via $g$) and applying $H_{\mathfrak{r}}^N(-)$ gives the following $B$-linear map:
\[
	H_{\mathfrak{r}}^N(g): \ H_{\mathfrak{r}}^N(B) \to H_{\mathfrak{r}}^N(B^g) = \ H_{\mathfrak{r}}^N(B)^g.
\]	
Now, Lemma \ref{lemma63} translates this $g$-linear map from $H_{\mathfrak{r}}^N(B) \to  H_{\mathfrak{r}}^N(B)$ into
\[
	H_{\mathfrak{r}}^N(g): B^\vee_B \to \ B^\vee_B.
\]
The previous Remark \ref{remark64} gives us $\lambda_i$ for $i = 1,\ldots,n$ such that $H_{\mathfrak{r}}^N(g)(\beta^\vee) = \lambda_s (\beta^\vee \circ g^{-1})$ for $\beta \in Be_s$.
\end{proof}
\end{lemma}
Finally, we can put the pieces together. Every graded automorphism $g$ of $B$ fixing the idempotents induces a $g$-linear map $H_{\mathfrak{r}}^N(B) \to H_{\mathfrak{r}}^N(B)$. Therefore, it comes along with $n$ scalars related to $H_{\mathfrak{r}}^N(g)(e_i^\vee)$. Putting the inverses of these scalars into a diagonal matrix finishes the definition of the homological determinant:
\begin{definition}
Let $g \in \Aut_{\gr}(B)$ as before. Define the \textit{homological determinant} $\Hdet(g)$ to be the diagonal matrix with entries $(\Hdet(g))_{ii} = \lambda_i^{-1}$ for the $\lambda_i$'s from Lemma \ref{lemma65}.
\end{definition}

\begin{proposition}
$\Hdet(-)$ forms a group homomorphism, i.e.
\[
	\Hdet(gh) = \Hdet(g)\Hdet(h)
\]
for all $g, h \in \Aut_{gr}(B)$.

\begin{proof}
Since both $g$ and $h$ fix the primitive orthogonal idempotents it is enough to consider their restrictions $g_i$ and $h_i$ and focus on $\Hdet(g)_{ii}$ and $\Hdet(h)_{ii}$ for some $1 \leq i \leq n$. Now, a similar proof to \cite[Proposition 2.5]{MR1758250} works. We obtain the following commutative triangle of $B$-linear maps
\[
\begin{tikzcd}[row sep = 30pt, column sep = 40pt]
e_iB \arrow[r, "h_i"] \arrow[rd, swap, "(gh)_i = g_ih_i"] & \left(e_iB\right)^{h} \arrow[d, "g_i"] \\
& \left(e_iB\right)^{gh} 
\end{tikzcd}
\]
Applying the functor $H_{\mathfrak{r}}^N$ gives
\[
\begin{tikzcd}[row sep = 40pt, column sep = 40pt]
H_{\mathfrak{r}}^N(e_iB) \arrow[r, "H_{\mathfrak{r}}^N(h_i)"] \arrow[rd, swap, "H_{\mathfrak{r}}^N(g_ih_i)"] & H_{\mathfrak{r}}^N\left((e_iB)^{h}\right) \arrow[d, "H_{\mathfrak{r}}^N(g_i)"] \\
& H_{\mathfrak{r}}^N\left((e_iB)^{gh}\right) 
\end{tikzcd}
\]
Lemma \ref{lemma65} and the definition of the homological determinant state that 
\[
	H_{\mathfrak{r}}^N(g_i)(-) = \Hdet(g)_{ii}^{-1} (g_i^{-1})^\ast(-)
\]
which equals scalar multiplication by $\Hdet(g)_{ii}^{-1}$ on the highest non-vanishing degree $-\ell_i$. Analogous observations are true for $H_{\mathfrak{r}}^N(h_i)(-)$ and $H_{\mathfrak{r}}^N(g_ih_i)(-)$. But the diagram commutes and therefore
\[
	\Hdet(g)_{ii}^{-1}\Hdet(h)_{ii}^{-1} = \Hdet(gh)_{ii}^{-1}.
\]
Since $\Hdet(-)$ equals the diagonal matrix with entries $\Hdet(-)_{ii}$ the claim follows.
\end{proof}
\end{proposition}

From the definition it is unclear how to calculate the $\lambda_i$'s or whether they are related to each other. In Section \ref{section7} we will prove that for preprojective algebras of quivers of extended Dynkin type $A, D$ or $E$ and graded automorphisms which scale the arrows, the homological determinant is a scalar matrix. However, this does not always need to be the case as one can see by looking at a direct product of two algebras.\\ \indent 
While we do not need to worry about shifts in the definition of the homological determinant, we now need to be more careful again. The following lemma follows \cite[Lemma 3.1]{MR1758250} closely and lays the groundwork for showing that $B^G$ is generalized Gorenstein when every $g \in G$ satisfies $\Hdet(g) = I_n$.

\begin{lemma}\label{lemma67}	Let $G$ be a finite group acting on $B$ such that every element fixes the primitive idempotents $e_i$ for $i = 1,\ldots,n$. We denote the fixed ring $B^G$ by $C$ and call $\mathfrak{r} = B_{\geq 1}$ and $\mathfrak{s} = C_{\geq 1}$.
\begin{itemize}[topsep=0pt]
	\item[(1)]	Let $N = \injdim_B(B)$. Then $H_{\mathfrak{s}}^i(C) = H_{\mathfrak{s}^{\op}}^i(C) = 0$ for $i \neq N$. In other words, $C$ is AS-Cohen-Macaulay.
	\item[(2)]	Assume every $g \in G$ satisfies $\Hdet(g) = I_n$. Then $H_\mathfrak{s}^N(C) \cong \bigoplus\limits_{b=1}^{n}{e_bC^\vee[\ell_b]}$ as right $C$-modules where $e_bC^\vee = \GrHom_k(Ce_b,k)$.
\end{itemize}
\begin{proof}
We start with (1). Recall that both $B$ and $C$ are noetherian by Lemma \ref{lemma211}(b) and in this section we assume that both $B$ and $B^{\op}$ satisfy the $\chi$-condition. This allows us to use results from \cite{MR1304753}. In particular, \cite[Lemma 8.2]{MR1304753} including its preparatory work holds in our case. Consequently, \cite[Lemma 4.15]{MR1674648} (Note: In \cite{MR1674648}, Lemma 4.15 is labeled 4.5 and appears between Theorem 4.14 and Theorem 4.16) is satisfied under our assumptions. This implies $H_{\mathfrak{r}}^i(M) = H_{\mathfrak{s}}^i(M)$ and $H_{\mathfrak{r}^{\op}}^i(N) = H_{\mathfrak{s}^{\op}}^i(N)$ for all graded right $C$-modules $M$ and graded left $C$-modules $N$ and all $i$. Thanks to Lemma \ref{lemma211}(a) we have $B \cong C \oplus K$ as $(C,C)$-bimodules and so
		\[
			H_{\mathfrak{s}}^i(C) \oplus H_{\mathfrak{s}}^i(K) = H_{\mathfrak{s}}^i(C \oplus K) = H_{\mathfrak{s}}^i(B) = H_{\mathfrak{r}}^i(B) = \begin{cases} 0 & \text{ for } i \neq N, \\ \bigoplus_{b=1}^{n}{e_bB^\vee}[\ell_b]  & \text{ for } i = N \end{cases} 
		\]
		by Lemma \ref{lemma63}. This shows $H_{\mathfrak{s}}^i(C) = 0$ for $i \neq N$, in other words, $C$ is AS-Cohen-Macaulay. The same calculation works over the opposite rings which finishes the first part.\\ \indent
Now, we prove (2). Recall from Lemma \ref{lemma211}(a) that the projection $\pi_G: B \to B^G = C$ via $b \mapsto \frac{1}{|G|}{\sum\limits_{g\in G}{g(b)}}$ is a splitting of the natural inclusion of $C$ into $B$, so we have 
\[
\begin{tikzcd}
0 \arrow[r] & C \arrow[rr,"\iota"] & & B \arrow[ll,"\pi_G",shift left = 1ex] \arrow[r] & B/C \arrow[r] & 0.
\end{tikzcd}
\]	
Being an additive functor,
		\[
			H_{\mathfrak{s}}^N(\pi_G): \ H_{\mathfrak{s}}^N(B) \to H_{\mathfrak{s}}^N(C)
		\]
		is the splitting of the inclusion of the right $C$-modules $H_{\mathfrak{s}}^N(C)$ into $H_{\mathfrak{s}}^N(B)$:
\[
\begin{tikzcd}
0 \arrow[r] & H_{\mathfrak{s}}^N(C) \arrow[rr,"\iota'"] & & H_{\mathfrak{s}}^N(B) \arrow[ll,"H_{\mathfrak{s}}^N(\pi_G)",shift left = 1ex] \arrow[r] & H_{\mathfrak{s}}^N\left(B/C\right) \arrow[r] & 0.
\end{tikzcd}
\]		
Then $H_{\mathfrak{s}}^N(\pi_G)\circ \iota' = \id_{H_{\mathfrak{s}}^N(C)}$ allows us to understand the module $H_{\mathfrak{s}}^N(C)$ as the fixed points of $H_{\mathfrak{s}}^N(\pi_G)$ via $\iota'$, i.e.
		\[
			H_{\mathfrak{s}}^N(C) = \{z \in H_{\mathfrak{s}}^N(B) \mid H_{\mathfrak{s}}^N(\pi_G)z = z\}.
		\]
		Next, we are going to show the equality
		\[
			H_{\mathfrak{s}}^N(C) = \left\{z \in H_{\mathfrak{s}}^N(B) \ \left| \ H_{\mathfrak{s}}^N(g)(z) = z \text{ for all } g \in G \right.\right\}
		\]
	using the identity
		\[
			H_{\mathfrak{s}}^N(\pi_G) = \frac{1}{|G|}\sum_{g \in G}{H_{\mathfrak{s}}^N(g)}.
		\]
	\begin{itemize}[topsep=0pt]
		\item[``$\supseteq$"]	Pick $z \in H_{\mathfrak{s}}^N(B)$ which is fixed under all $H_{\mathfrak{s}}^N(g)$. Clearly, $z$ is also fixed under the weighted sum $H_{\mathfrak{s}}^N(\pi_G)$ and thus $z \in H_{\mathfrak{s}}^N(C)$.
		\item[``$\subseteq$"]	Conversely, let $z \in H_{\mathfrak{s}}^N(C)$ which says $H_{\mathfrak{s}}^N(\pi_G)(z) = z$. The same calculation as in \cite[Lemma 3.1]{MR1758250} provides
		\begin{align*}
			H_{\mathfrak{s}}^N(g)(z) &= H_{\mathfrak{s}}^N(g)\left(H_{\mathfrak{s}}^N(\pi_G)(z)\right) = H_{\mathfrak{s}}^N(g) \sum_{h \in G}{\frac{1}{|G|}H_{\mathfrak{s}}^N(h)(z)} \\ 
			&= \sum_{h \in G}{\frac{1}{|G|}H_{\mathfrak{s}}^N(gh)(z)} = H_{\mathfrak{s}}^N(\pi_G)(z) = z
		\end{align*}
		for all $g$. This proves the second inclusion.
	\end{itemize}

	\indent	Combined, together with the observation $H_{\mathfrak{r}}^N(M) = H_{\mathfrak{s}}^N(M)$ for all graded right $C$-modules $M$ made in part (1), it follows that
\begin{equation}
	\begin{aligned}\label{eq3-8}
		H_{\mathfrak{s}}^N(C) &= \left\{z \in H_{\mathfrak{r}}^N(B) \ \left| \ H_{\mathfrak{s}}^N(g)(z) = z \text{ for all } g \in G \right. \right\} \\
		&\overset{\ref{lemma63}}{=} \left\{ \left. z \in \bigoplus_{b=1}^{n}{e_bB^\vee}[\ell_b] \ \right| \ H_{\mathfrak{s}}^N(g)(z) = z \text{ for all } g \in G \right\}.
	\end{aligned}
\end{equation}
The assumption $\Hdet(g) = I_n$ implies that $H_{\mathfrak{s}}^N(g)(\beta^\vee) = (\beta^\vee \circ g^{-1}) = (g^{-1})^\ast(\beta^\vee)$ for all $g \in G$ and all $\beta \in B$. We apply this to Equation \eqref{eq3-8}:
	\[
	H_{\mathfrak{s}}^N(C) = \left\{ \left. z \in \bigoplus_{b=1}^{n}{e_bB^\vee}[\ell_b] \ \right| \ (g^{-1})^\ast(z) = z \text{ for all } g \in G \right\} \overset{(\ast)}{=} \bigoplus\limits_{b=1}^{n}{e_bC^\vee[\ell_b]}.
	\]
	To see $(\ast)$, recall that any element in $C$ is fixed under the action of $G$. Therefore, every $\beta^\vee \in \ e_bC^\vee$ satisfies $(g^{-1})^\ast(\beta^\vee) = \beta^\vee \circ g^{-1} = \beta^\vee$. The reverse argument is similar. 
\end{proof}
\end{lemma}

\begin{lemma}\label{lemma68}
Let $C$ be a noetherian $\NN$-graded AS-Cohen-Macaulay algebra with $C_0 \cong k^n$ and denote its primitive pairwise orthogonal idempotents by $\{e_1,\ldots,e_n\}$. Assume that $C$ has a balanced dualizing complex. We write $\mathfrak{s} = C_{\geq 1}$. If $C$ satisfies $H_{\mathfrak{s}}^N(C) = \bigoplus\limits_{b=1}^{n}{e_bC^\vee[\ell_b]}$ as right $C$-modules, then $C$ is generalized Gorenstein.
\begin{proof}
Abusing notation, we denote the simple right and left $C$-modules by $S_i = e_iC/(e_iC)_{\geq 1}$ and $S_i^\vee = Ce_i/(Ce_i)_{\geq 1}$, respectively. 
Similarly to \cite[Lemma 3.2]{MR1758250}, the modified spectral sequence from \cite[Proposition 1.1]{MR1643863} given by
\[
	E_2^{pq} = \GrExt_C^p(C_0,H_\mathfrak{s}^q(C)) \ \R \ \GrExt_C^{p+q}(C_0,C)
\]
collapses due to $C$ being AS-Cohen-Macaulay. Moreover, as used before, all the $e_bC^\vee$ are injective and so $E^{pq}$ has only one nonzero entry at the $p = 0$ and $q = N$ position. We are able to draw two conclusions. First, $\GrExt_C^i(C_0,C) = 0$ for $i \neq N$ and second,
\begin{equation}\label{eq5}
	\GrExt_C^N(C_0,C) \cong \GrExt_C^0(C_0,H_\mathfrak{s}^N(C)) = \GrHom_C\left(C_0, \ \bigoplus\limits_{b=1}^{n}{e_bC^\vee[\ell_b]}\right) = \bigoplus_{b=1}^{n}{S_b^\vee[\ell_b]}.
\end{equation}
By a symmetric version of \cite[Lemma 2.8]{DanToDo} it follows that
\[
	\GrExt_C^N(e_iC_0,C) \cong \bigoplus_{b=1}^{n}{S_b^\vee[\ell_b]}e_i \neq \{0\}
\]
for all $i$ and together with Equation \eqref{eq5} this implies $\dim_k(\GrExt_C^N(S_i,C)) = 1$. In other words, $C$ satisfies the generalized Gorenstein condition.\\ \indent
To finish the proof, we need to show that $\injdim_C(C)$ is finite. For this, we apply \cite[Theorem 6.3]{MR1469646} which works for non-connected algebras as explained in the first paragraph of \cite[Lemma 3.5]{MR3250287}. We keep in mind that $C$ has a balanced dualizing complex and so the aforementioned result verifies both requirements of van den Bergh's Local Duality Theorem \cite[Theorem 5.1]{MR1469646}. In particular, the $\chi$-condition implies that $\Ext^i(C_0,C_0)$ is finite dimensional for all $i$ which is the non-connected equivalent of being $\Ext$-finite which is one of the requirements. Finally, $H_{\mathfrak{s}}^N(C)^\vee$ has finite injective dimension by \cite[Theorem 5.1(1)]{MR1469646} and $H_{\mathfrak{s}}^N(C)^\vee = C^{\vee\vee} = C$ finishes the proof.
\end{proof}
\end{lemma}
At last, we can put together the previous lemmas. Fortunately, fixed rings have a balanced dualizing complex and so the following main theorem is the combination of \ref{lemma67} and \ref{lemma68}.
\begin{theorem}\label{thm69}
Suppose $B$ has finite global dimension $d$. Suppose $G$ is a finite subgroup of graded automorphisms of $B$ such that every $g \in G$ fixes the primitive idempotents. Assume that every $g \in G$ satisfies $\Hdet(g) = I_n$. Then $B^G$ is generalized Gorenstein.
\begin{proof}
We use the notation $C = B^G$, $\mathfrak{s} = C_{\geq 1}$ and $\mathfrak{r} = B_{\geq 1}$ as before. By Lemma \ref{lemma67} it follows that $C$ is AS-Cohen-Macaulay and that $H_{\mathfrak{s}}^N(C) \cong \bigoplus\limits_{b=1}^{n}{e_bC^\vee[\ell_b]}$. Thanks to Proposition \ref{prop213} the algebra $C$ satisfies the $\chi$-condition. Our assumptions and Lemma \ref{lemma211}(b) guarantee that both $B$ and $C$ are noetherian.\\ \indent
We want to again apply \cite[Theorem 6.3]{MR1469646} to $B$ which holds for non-connected algebras. The first condition is satisfied since $B$ has finite global dimension. This instantly implies that the functors 
\begin{align*}
&H_\mathfrak{r}^0(M): \GrMod\text{-}B \longrightarrow \GrMod\text{-}B, \ \ \ M \mapsto \varinjlim \Hom_B(B/B_{\geq r}, M),
\end{align*}
and similarly $\left(H_\mathfrak{r}^0(M)\right)^{\op}$ for left $B$-modules have finite cohomological dimension. Moreover, $B$ and $B^{\op}$ satisfy the $\chi$-condition by hypothesis. Hence, $B$ has a balanced dualizing complex. Using the projection $B \to B^G$ which is a finite homomorphism of graded $k$-algebras, \cite[Theorem 4.16]{MR1674648} yields a balanced dualizing complex for $B^G$. With the requirements of Lemma \ref{lemma68} satisfied, $C$ is generalized Gorenstein.
\end{proof}
\end{theorem}
Notice that $B$ and $B^{\op}$ in Theorem \ref{thm69} do not need to be assumed in advance to satisfy the $\chi$-condition. Since $B$ has finite global dimension and satisfies the generalized Gorenstein condition, it must automatically satisfy $\chi$ due to \cite[Lemma (2.8)]{OtherPaper}.

\section{The vector trace of an automorphism of a preprojective algebra}\label{sectionVT}

\subsection{Case $\widetilde{A_{n-1}}$}

Let $n \geq 2$ and let $Q$ be the following quiver coming from the extended Dynkin diagram $\widetilde{A_{n-1}}$:
\[
\begin{tikzcd}
& & & n \ \bullet \arrow[llldd,"\alpha_{n}"] \\
& \\
1 \ \bullet \arrow[r,swap,"\alpha_1"] & \bullet \ 2 \arrow[r,swap,"\alpha_2"] & \bullet \ 3 & \ldots & \bullet \arrow[r] & \bullet \arrow[r] &  \arrow[llluu, "\alpha_{n-1}"] \bullet \ n-1
\end{tikzcd}
\]

\vspace{12pt}
In this subsection, let $A = \Pi(Q)$ be the preprojective algebra with respect to $Q$. Mainly, we are stating formulas for $\TrV(g,t)$ and $\TrVL(g,t)$ for any graded automorphism $g$ of $A$ which satisfies Convention \ref{con214} in Proposition \ref{prop33}. In this situation, we can denote $g(\alpha_i) = c_i \alpha_i$ and $g(\alpha_i^\ast) = t_i \alpha_i^\ast$. 

\vspace{12pt}
\begin{defremark}\label{remark31}	
It is easy to see that $A$ can be understood as the path algebra $k\bar{Q}$ with respect to the relations $\alpha_i \alpha_i^\ast = \alpha_{i-1}^\ast \alpha_{i-1}$ for $i = 1,\ldots,n$ and reduce the indices modulo $n$ if needed.  \\ \indent
As a consequence, we can always rearrange a path in $A$ to first have the nonstar arrows occur followed by the star arrows. Counting as done in \cite[(3.1) Definition \& Remark]{OtherPaper} yields:
\[
	\HS_{e_\ell A}(t) = 1 + 2t + 3t^2 + 4t^3 + \ldots = \sum_{j = 0}^{\infty}{(j+1)t^j} = \frac{1}{(1-t)^2}. 
\]
Every element of $A$ can be represented as a $k$-linear combination of paths in the double of $Q$. In order to simplify notation, we identify an element of $A$ with any representative in $k\bar{Q}$. We call an element of the form $\prod_{j}{\alpha_{i_j}}\prod_{v}{\alpha_{u_v}^\ast}$ a \textit{simple path}. Due to the relations, every path has a unique representation as a simple path. Moreover, every simple path in $A$ is uniquely determined by its starting vertex, its length and the number of star arrows. Further, every element of $A$ is a unique linear combination of simple paths.
\end{defremark}

The following proposition from {\cite{dissertation}} presents formulas for the vector traces of a graded automorphism of $A$. These will be used in Section \ref{section7} to find a reciprocity statement connecting $\TrV(g,t)$ to $\TrVL(g^{-1},t^{-1})$ which eventually helps to find concrete formulas for $\Hdet(g)$. 

\begin{proposition}[{\cite[Proposition 4.1.5]{dissertation}}]\label{prop33} Let $Q = \widetilde{A_{n-1}}$ and $A = \Pi(Q)$. Fix $g$ to be a graded automorphism of $A$ with $g(\alpha_i) = c_i \alpha_i$ and $g(\alpha_i^\ast) = t_i \alpha_i^\ast$. Then for $n \geq 3$
	\[
	\TrV(g,t) = \underbrace{\begin{pmatrix} 1+c_1t_1t^2 & -c_1t &  & & & -t_nt \\[2ex] -t_1t & 1 + c_1t_1t^2 & -c_2t & & & \\[2ex] & -t_2t & 1+c_1t_1t^2 & -c_3 t & & \\[2ex] & & \ddots & \ddots & \ddots & \\[2ex]  & & & \ddots & \ddots & -c_{n-1} t \\[2ex] -c_nt & & & & -t_{n-1}t & 1+c_1t_1t^2 \end{pmatrix}^{-1}}_{= M(c_i ,t_i ,t)^{-1}} \cdot \begin{pmatrix} 1 \\[2ex] 1 \\[2ex] 1 \\[2ex] \vdots \\[2ex] 1 \\[2ex] 1 \end{pmatrix}
	\]
	and
	\[
	\TrVL(g,t) = M(t_i,c_i,t)^{-1} \cdot \begin{pmatrix} 1 \\ \vdots \\ 1 \end{pmatrix}
	\]
	where $M(c_i,t_i,t)^T = M(t_i,c_i,t)$. In the special case $n = 2$, which can be easily extracted from the proof of \cite[Proposition 4.1.5]{dissertation}, we have instead
		\[
	\TrV(g,t) = \underbrace{\begin{pmatrix} 1+c_1t_1t^2 & -(c_1+t_2)t \\[2ex] -(c_2+t_1)t & 1+c_1t_1t^2 \end{pmatrix}^{-1}}_{= M(c_i ,t_i ,t)^{-1}} \cdot \begin{pmatrix} 1 \\[2ex] 1 \end{pmatrix}
	\]
	and
	\[
	\TrVL(g,t) = M(t_i,c_i,t)^{-1} \cdot \begin{pmatrix} 1 \\ 1 \end{pmatrix} = M(c_i,t_i,t)^{-T} \cdot \begin{pmatrix} 1 \\ 1 \end{pmatrix}.
	\]
\end{proposition}

\subsection{Case $\widetilde{D_{n-1}}$}
Next, let $n \geq 5$ and let $Q$ be the following quiver coming from the extended Dynkin diagram $\widetilde{D_{n-1}}$:
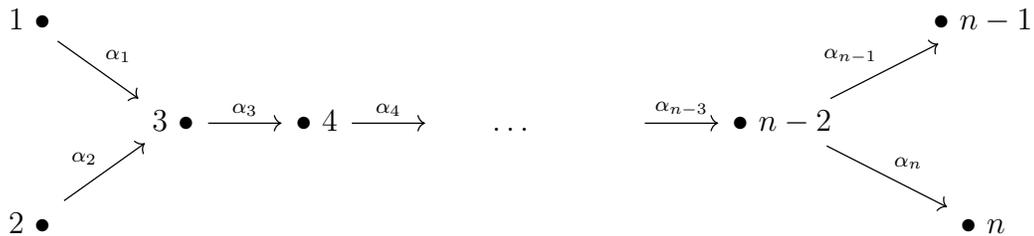
\[
\begin{tikzcd}
1 \ \bullet \arrow[rd,"\alpha_1"]& & & & & & \bullet \ n-1 \\
& 3 \ \bullet \arrow[r,"\alpha_3"] & \bullet \ 4 \arrow[r,"\alpha_4"] & \hspace{20pt} \ldots & \arrow[r,"\alpha_{n-3}"] & \bullet \ n-2  \arrow[ru,"\alpha_{n-1}"] \arrow[rd,"\alpha_n"] \\
2 \ \bullet \arrow[ru,"\alpha_2"] & & & & & & \bullet \ n
\end{tikzcd}
\]
As before, let $A = \Pi(Q)$ be the preprojective algebra with respect to $Q$. We need to give a description of $\TrV(g,t)$ and $\TrVL(g,t)$ of a graded automorphism $g$ as in Convention \ref{con214} in this case. But first, we analyze the relations of $A$ as they are crucial to find the minimal projective resolutions of the simple modules.

\begin{remark}\label{remark41}	In Definition \ref{def23}, the preprojective algebra was defined as the path algebra modulo the ideal $I$ generated by the one relation:
\begin{align}\label{eqrelD}
\sum_{\alpha \in Q_1}{\alpha \alpha^\ast - \alpha^\ast \alpha}.
\end{align}
By successively multiplying Equation \eqref{eqrelD} by the trivial paths $e_i$ for $i = 1,\ldots,n$ we obtain the relations
\begin{align}\label{eqrelDalle}
	\alpha_1 \alpha_1^\ast &= 0, \notag \\
	\alpha_2 \alpha_2^\ast &= 0, \notag \\
	\alpha_3 \alpha_3^\ast &= \alpha_1^\ast \alpha_1 + \alpha_2^\ast \alpha_2, \notag \\
	\alpha_i \alpha_i^\ast &= \alpha_{i-1}^\ast \alpha_{i-1}, \ \ \ \ \ \ \ \ \ \ \ \ \ \ \ \ \ i = 4,\ldots,n-3, \\
	\alpha_{n-3}^\ast \alpha_{n-3} &= \alpha_{n-1}\alpha_{n-1}^\ast + \alpha_n \alpha_n^\ast, \notag \\
	\alpha_{n-1}^\ast \alpha_{n-1} &= 0,	\hspace{4cm} \text{ and } \notag \\
	\alpha_n^\ast \alpha_n &= 0. \notag
\end{align}
Starting now we think of $I$ being generated by the relations in Equation \eqref{eqrelDalle}. Moreover, since $g$ scales the arrows, notice that these relations force $c_it_i = c_1t_1$ for all $i = 1,\ldots,n$ as it was the case for $Q = \widetilde{A_{n-1}}$. Conversely, any nonzero scalars $c_i$, $t_i$ satisfying
$c_it_i = c_1t_1$ give rise to a graded automorphism of $A$.
\end{remark}
For the proof of the subsequent proposition, we follow \cite[Theorem 2.3]{MR1438180} closely.

\begin{proposition}\label{prop42}
Let $Q = \widetilde{D_{n-1}}$ for $n \geq 5$ and $A = \Pi(Q)$. Fix $g$ to be a graded automorphism of $A$ with $g(\alpha_i) = c_i \alpha_i$ and $g(\alpha_i^\ast) = t_i \alpha_i^\ast$ for $1 \leq i \leq n-3$ and $i = n-1, n$. Then, $\TrV(g,t)$ equals
\[
	\underbrace{\begin{pmatrix} 1+c_1t_1t^2 & 0 & -c_1t & & & & \\[1ex] 0 & 1 + c_1t_1t^2 & -c_2t & & & & \\[1ex] -t_1t & -t_2t & 1+c_1t_1t^2 & -c_3 t & & & \\[1.5ex] & & \ddots & \ddots & \ddots & & \\[1.5ex] & & & -t_{n-3}t & 1+c_1t_1t^2 & -c_{n-1}t & -c_nt \\[1.5ex]  & & & & -t_{n-1}t & 1+c_1t_1t^2 & 0 \\[1ex] & & & & -t_nt & 0 & 1+c_1t_1t^2 \end{pmatrix}^{-1}}_{= M(c_i,t_i,t)^{-1}} \cdot \begin{pmatrix} 1 \\[1.5ex] 1 \\[2ex] \vdots \\[2ex] 1 \\[1.5ex] 1 \end{pmatrix}
\]
and similarly 
\[ 
\TrVL(g,t) = M(t_i,c_i,t)^{-1} \cdot \begin{pmatrix} 1 \\ \vdots  \\ 1 \end{pmatrix}
\]
where $M(c_i,t_i,t)^T = M(t_i,c_i,t)$.
\begin{proof}
We wish to apply \cite[Lemma 2.1]{MR1438180} to determine the trace of $g$. To this end, we first write down a projective resolution of each $S_i$. Since $g$ preserves the indecomposable projective modules $P_i = e_iA$ as well as the graded submodules $e_iA_{\geq 1}$, one can lift the minimal projective resolutions of every simple module $S_i =  e_iA / e_iA_{\geq 1}$ to a minimal projective resolution of $S_i^g$. Then we can use that the alternating sum of the connecting homomorphisms vanishes and deduce a formula of the trace functions. \\ \indent
Due to $\gldim(A) = 2$, \cite[Proposition 4.2]{MR1930968} gives the required minimal projective resolution of $S_i$ for $i = 4,\ldots,n-3$:
\[
	0 \longrightarrow e_iA[-2] \overset{\delta_2 = \begin{pmatrix} \alpha_i^\ast \\ -\alpha_{i-1}\end{pmatrix} \boldsymbol{\cdot}}{\xrightarrow{\hspace*{3cm}}} e_{i+1}A[-1] \oplus e_{i-1}A[-1] \overset{\delta_1 = \begin{pmatrix} \alpha_i & \alpha_{i-1}^{\ast} \end{pmatrix} \boldsymbol{\cdot}}{\xrightarrow{\hspace*{3cm}}} e_iA \overset{\delta_0}{\longrightarrow} S_i \longrightarrow 0.
\]
We define $g_i = g|_{e_iA}$ and let $\tilde{g}_i$ be the induced map of $g_i$ on $S_i$. Then $\tilde{g}_i$ is a $g$-linear map from $S_i$ to $S_i$ or an $A$-linear map from $S_i$ to $S_i^g$. We consider the following:
\[
\begin{tikzcd}
0 \arrow[r] & Q_2 = e_iA \arrow[d, "\phi_2^{(i)}"] \arrow[r, "\delta_2"] & Q_1 = e_{i+1}A \oplus e_{i-1}A \arrow[d, "\phi_1^{(i)}"]\arrow[r, "\delta_1"] & Q_0 = e_iA \arrow[d, "\phi_0^{(i)}"] \arrow[r, "\delta_0"] & S_i \arrow[r] \arrow[d, "\tilde{g}_i"] & 0 \\
0 \arrow[r] & Q_2^g = (e_iA)^g \arrow[r, "\delta_2"] & Q_1^g = (e_{i+1}A \oplus e_{i-1}A)^g \arrow[r, "\delta_1"] & Q_0^g = (e_iA)^g \arrow[r, "\delta_0"] & S_i^g \arrow[r] & 0.
\end{tikzcd}
\]
Notice that the new multiplicative structure on $Q_j^g$ takes part in the action of $\delta_j: Q_j^g \to Q_{j-1}^g$. For instance, $\delta_1: Q_1^g \to Q_0^g$ acts via 
\[
	\delta_1 \begin{pmatrix} e_{i+1} \ast a \\ e_{i-1} \ast b \end{pmatrix} = e_i \alpha_i g(b) + e_{i} \alpha_{i-1}^\ast g(b).
\]
Using the universal property of projective resolutions gives graded $A$-linear maps $\phi_j^{(i)}: Q_j \to Q_j^g$ such that the squares commute.\\ \indent
Recall the relations of $A$ given in Equation \eqref{eqrelDalle}. Simple calculations give $\phi_0^{(i)} = g_i$, $\phi_1^{(i)}\begin{pmatrix} e_{i+1} \\ 0 \end{pmatrix} = \begin{pmatrix} c_ie_{i+1} \\ 0 \end{pmatrix}$ and $\phi_1^{(i)}\begin{pmatrix} 0 \\ e_{i-1} \end{pmatrix} = \begin{pmatrix} 0 \\ t_{i-1} e_{i-1} \end{pmatrix}$ as well as $\phi_2^{(i)}(e_i) = c_1t_1 e_i$.\\ \indent
Using the same technique, we obtain the following commutative diagrams whose rows are projective resolutions
\[
\begin{tikzcd}
0 \arrow[r] & Q_2 = e_jA \arrow[d, "\phi_2^{(j)}"] \arrow[r, "\alpha_j^\ast \boldsymbol{\cdot}"] & Q_1 = e_3A \arrow[d, "\phi_1^{(j)}"]\arrow[r, "\alpha_j \boldsymbol{\cdot}"] & Q_0 = e_j A \arrow[d, "\phi_0^{(j)}"] \arrow[r, "\delta_0"] & S_j \arrow[r] \arrow[d, "\tilde{g}_j"] & 0 \\
0 \arrow[r] & Q_2^g = (e_jA)^g \arrow[r, "\alpha_j^\ast \boldsymbol{\cdot}"] & Q_1^g = (e_3A)^g \arrow[r, "\alpha_j \boldsymbol{\cdot}"] & Q_0^g = (e_jA)^g \arrow[r, "\delta_0"] & S_j^g \arrow[r] & 0
\end{tikzcd}
\]
for $j = 1,2$ with $\phi_0^{(j)} = g_j, \ \phi_1^{(j)}(e_3) = c_j e_3$ as well as $\phi_2^{(j)}(e_j) = c_1t_1 e_j$. Similarly,
\[
\begin{tikzcd}
0 \arrow[r] & Q_2 = e_lA \arrow[d, "\phi_2^{(l)}"] \arrow[r, "\alpha_l \boldsymbol{\cdot}"] & Q_1 = e_{n-2}A \arrow[d, "\phi_1^{(l)}"]\arrow[r, "\alpha_l^\ast \boldsymbol{\cdot}"] & Q_0 = e_l A \arrow[d, "\phi_0^{(l)}"] \arrow[r, "\delta_0"] & S_l \arrow[r] \arrow[d, "\tilde{g}_l"] & 0 \\
0 \arrow[r] & Q_2^g = (e_lA)^g \arrow[r, "\alpha_l \boldsymbol{\cdot}"] & Q_1^g = (e_{n-2}A)^g \arrow[r, "\alpha_l^\ast \boldsymbol{\cdot}"] & Q_0^g = (e_lA)^g \arrow[r, "\delta_0"] & S_l^g \arrow[r] & 0
\end{tikzcd}
\]
for $l = n-1,n$ with $\phi_0^{(l)} = g_l, \ \phi_1^{(l)}(e_{n-2}) = t_l e_{n-2}$ as well as $\phi_2^{(l)}(e_l) = c_1t_1 e_l$ are obtained.\\ \indent
The two remaining resolutions are a little more complicated but can also be computed using that the squares commute and that all maps are $g$-linear from $Q_r$ to $Q_r$. First we have
\begin{small}
\[
\begin{tikzcd}
0 \arrow[r] & Q_2 = e_3A \arrow[d, "\phi_2^{(3)}"] \arrow[r, "\delta_2 = \begin{pmatrix} \alpha_3^\ast \\ -\alpha_{2} \\ -\alpha_1 \end{pmatrix} \boldsymbol{\cdot}"] & Q_1 = e_4A \oplus e_2A \oplus e_1A \arrow[d, "\phi_1^{(3)}"]\arrow[r, "\delta_1 = \begin{pmatrix} \alpha_3 \\ \alpha_2^{\ast} \\ \alpha_1^{\ast} \end{pmatrix}^T \boldsymbol{\cdot}"] & Q_0 = e_3A \arrow[d, "\phi_0^{(3)}"] \arrow[r, "\delta_0"] & S_3 \arrow[r] \arrow[d, "\tilde{g}_3"] & 0 \\
0 \arrow[r] & Q_2^g = (e_3A)^g \arrow[r, "\delta_2"] & Q_1^g = (e_4A \oplus e_2A \oplus e_1A)^g \arrow[r, "\delta_1"] & Q_0^g = (e_3A)^g \arrow[r, "\delta_0"] & S_3^g \arrow[r] & 0
\end{tikzcd}
\]
\end{small}
with $\phi_0^{(3)} = g_3$, $\phi_2^{(3)}(e_3) = c_1t_1 e_3$ and
\[
	\phi_1^{(3)}\begin{pmatrix} e_4 \\ 0 \\ 0 \end{pmatrix} = \begin{pmatrix} c_3e_4 \\ 0 \\ 0 \end{pmatrix}, \ \ \ \ \ \phi_1^{(3)}\begin{pmatrix} 0 \\ e_2 \\ 0 \end{pmatrix} = \begin{pmatrix} 0 \\ t_2e_2 \\ 0 \end{pmatrix}, \ \ \ \ \ \ \phi_1^{(3)}\begin{pmatrix} 0 \\ 0 \\ e_1 \end{pmatrix} = \begin{pmatrix} 0 \\ 0 \\ t_1 e_1 \end{pmatrix}
\]
which uniquely determines the maps $\phi_0^{(3)}, \phi_1^{(3)}$ and $\phi_2^{(3)}$.\\ \indent
Finally, the last resolution computes in the same way yielding
\[
\begin{tikzcd}
0 \arrow[r] & e_{n-2}A \arrow[d, "\phi_2^{(n-2)}"] \arrow[r, "\delta_2 = \begin{pmatrix} \alpha_{n-1}^\ast \\ \alpha_{n}^\ast \\ -\alpha_{n-3} \end{pmatrix} \boldsymbol{\cdot}"] & e_{n-1}A \oplus e_nA \oplus e_{n-3}A \arrow[d, "\phi_1^{(n-2)}"]\arrow[r, "\delta_1 = \begin{pmatrix} \alpha_{n-1} \\ \alpha_n \\ \alpha_{n-3}^{\ast} \end{pmatrix}^T \boldsymbol{\cdot}"] & e_{n-2}A \arrow[d, "\phi_0^{(n-2)}"] \arrow[r, "\delta_0"] & S_{n-2} \arrow[r] \arrow[d, "\tilde{g}_{n-2}"] & 0 \\
0 \arrow[r] & (e_{n-2}A)^g \arrow[r, "\delta_2"] & = (e_{n-1}A \oplus e_nA \oplus e_{n-3}A)^g \arrow[r, "\delta_1"] & (e_{n-2}A)^g \arrow[r, "\delta_0"] & S_{n-2}^g \arrow[r] & 0
\end{tikzcd}
\]
with $\phi_0^{(n-2)} = g_{n-2}$, $\phi_2^{(n-2)}(e_{n-2}) = c_1t_1 e_{n-2}$ and
\begin{align*} 
	\phi_1^{(n-2)}\begin{pmatrix} e_{n-1} \\ 0 \\ 0 \end{pmatrix} = \begin{pmatrix} c_{n-1}e_{n-1} \\ 0 \\ 0 \end{pmatrix}, \
	\phi_1^{(n-2)}\begin{pmatrix} 0 \\ e_n \\ 0 \end{pmatrix} = \begin{pmatrix} 0 \\ c_ne_n \\ 0 \end{pmatrix} \text{ and } 
	\phi_1^{(n-2)}\begin{pmatrix} 0 \\ 0 \\ e_{n-3} \end{pmatrix} = \begin{pmatrix} 0 \\ 0 \\ t_{n-3} e_{n-3} \end{pmatrix}.
\end{align*}
Viewing the maps $\delta_i$ as graded homomorphisms induces a natural grading on the modules $Q_r$. That is, each $Q_r$ is the sum of shifted copies of projective modules $e_iA$. Further, the maps $\phi_r^{(s)}$ are graded maps and the mentioned shifts in the grading of the $Q_r$ and $Q_r^g$ become visible in the traces of the $\phi_r^{(s)}$. For $s = 1,\ldots,n$, we now view the $\phi_r^{(s)}$ as $g$-linear graded maps from $Q_r$ to $Q_r$. Since $\phi_0^{(s)} = g_s$, we find that $\Tr_{Q_0}(\phi_0^{(s)},t) = \Tr_{e_sA}(g,t)$. Applying $\phi_2^{(s)}$ to $\beta \in e_sA$ gives 
\[
\phi_2^{(s)}(\beta) = \phi_2^{(s)}(e_s)g(\beta) = c_1t_1 g(\beta)
\]
and therefore $\Tr_{Q_2}(\phi_2^{(s)},t) = \Tr_{e_sA[-2]}(c_1t_1g,t) = c_1t_1t^2 \Tr_{e_sA}(g,t).$ Notice that the $t^2$ comes from the shift in the grading of $Q_2$.\\ \indent
It remains to determine $\Tr_{Q_1}(\phi_1^{(s)},t)$. First, consider $s = 4,\ldots,n-3$. In the same way as before but using bases for each component individually, we immediately obtain the equation $\Tr_{Q_1}(\phi_1^{(s)},t) = c_it \Tr_{e_{i+1}A}(g,t) + t_{i-1}t \Tr_{e_{i-1}A}(g,t)$. From the projective resolutions above we read off the remaining four cases:
\begin{align}\label{eq3new}
\begin{split}
	\Tr_{Q_1}(\phi_1^{(j)},t) &= c_jt\Tr_{e_3A}(g,t), \ \ \ \ \ \ \text{ for } j = 1,2, \\	
	\Tr_{Q_1}(\phi_1^{(l)},t) &= t_lt \Tr_{e_{n-2}A}(g,t), \ \ \ \ \text{ for } l = n-1,n, \\
	\Tr_{Q_1}(\phi_1^{(3)},t) &= c_3t \Tr_{e_4A}(g,t) + t_2t \Tr_{e_2A}(g,t) + t_1t \Tr_{e_1A}(g,t), \ \ \ \text{ and } \\
	\Tr_{Q_1}(\phi_1^{n-2},t) &= c_{n-1}t \Tr_{e_{n-1}A}(g,t) + c_nt \Tr_{e_nA}(g,t) + t_{n-3}t\Tr_{e_{n-3}A}(g,t).
\end{split}
\end{align}
Now, \cite[Lemma 2.1]{MR1438180} applies and provides $\sum\limits_{j = 0}^{2}{(-1)^j \Tr_{Q_j}(\phi_j^{(i)},t)} = \Tr_{S_i}(\tilde{g}_i,t)$ which translates to:
\begin{align*}
	\Tr_{S_i}(\tilde{g}_i,t) &= \Tr_{e_iA}(g,t) - \left[c_it \Tr_{e_{i+1}A}(g,t) + t_{i-1}t \Tr_{e_{i-1}A}(g,t)\right] + c_1t_1t^2 \Tr_{e_iA}(g,t), \\ & \hspace{11cm} i = 4,\ldots,n-3, \\
	\Tr_{S_j}(\tilde{g}_j,t) &= \Tr_{e_jA}(g,t) - c_jt \Tr_{e_3A}(g,t) + c_1t_1t^2 \Tr_{e_jA}(g,t), \ \ \ \ \ \ j = 1,2, \\
	\Tr_{S_l}(\tilde{g}_l,t) &= \Tr_{e_lA}(g,t) - t_lt \Tr_{e_{n-2}A}(g,t) + c_1t_1t^2 \Tr_{e_lA}(g,t), \ \ \ \ l = n-1,n, \\
	\Tr_{S_3}(\tilde{g}_3,t) &= \Tr_{e_3A}(g,t) - \left[c_3t \Tr_{e_4A}(g,t) + t_2 t \Tr_{e_2A}(g,t) + t_1t \Tr_{e_1A}(g,t)\right] \\
	&\ \ \ + c_1t_1t^2 \Tr_{e_3A}(g,t), \ \ \ \text{ and } \\
	\Tr_{S_{n-3}}(\tilde{g}_{n-3},t) &= \Tr_{e_{n-3}A}(g,t) - \left[c_{n-1}t \Tr_{e_{n-1}A}(g,t) + c_nt \Tr_{e_nA}(g,t) + t_{n-3}t \Tr_{e_{n-3}A}(g,t)\right]\\
	&\ \ \ + c_1t_1t^2 \Tr_{e_{n-3}A}(g,t).
\end{align*}

Finally, we can write the $n$ equations using matrix notation as
\[
\begin{pmatrix} \Tr_{S_1}(\tilde{g}_1,t) & \Tr_{S_2}(\tilde{g}_2,t) & \ldots & \ldots & \ldots & \Tr_{S_n}(\tilde{g}_n,t) \end{pmatrix}^T = M(c_i,t_i,t) \boldsymbol{\cdot} \TrV(g,t),
\]
where
\[
M(c_i,t_i,t) = \begin{pmatrix} 1+c_1t_1t^2 & 0 & -c_1t & & & & \\[1ex] 0 & 1+c_1t_1t^2 & -c_2t & & & & \\[1ex] -t_1t & -t_2t & 1+c_1t_1t^2 & -c_3 t & & & \\[1.5ex] & & \ddots & \ddots & \ddots & & \\[1.5ex] & & & -t_{n-3}t & 1+c_1t_1t^2 & -c_{n-1}t & -c_nt \\[1.5ex]  & & & & -t_{n-1}t & 1+c_1t_1t^2 & 0 \\[1ex] & & & & -t_nt & 0 & 1+c_1t_1t^2 \end{pmatrix}.
\]
On the left hand side, $\Tr_{S_i}(\tilde{g}_i,t) = 1$ as $\tilde{g}_i$ sends $e_i$ to $e_i$. Therefore, multiplying the inverse of the matrix $M(c_i,t_i,t)$ from the left finishes the proof for $\TrV(g,t)$. The case $\TrVL(g,t)$ follows from analogous calculations.
\end{proof}
\end{proposition}

\subsection{Case $\widetilde{E_m}$ for $m = 6,7,8$}
Finally, the cases $\widetilde{E_m}$ for $m = 6,7,8$ remain. Let $g$ be as in Convention \ref{con214} for $A = \Pi(Q)$ where the underlying quivers $Q$ are as follows.
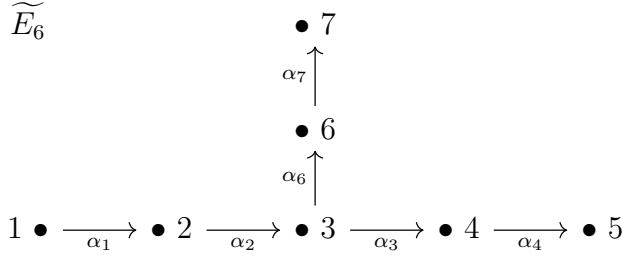
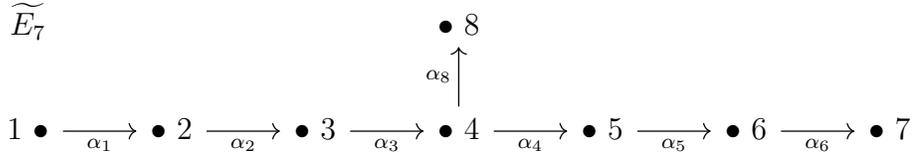
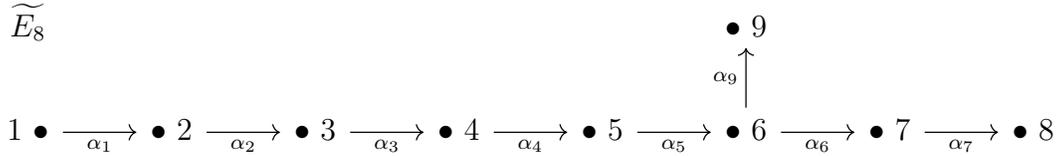
\[
\begin{tikzcd}
\widetilde{E_6} & & \bullet \ 7 \\
& & \bullet \arrow[u,"\alpha_7"] \ 6  \\
1 \ \bullet \arrow[r,swap,"\alpha_1"] & \bullet \ 2 \arrow[r,swap,"\alpha_2"] & \bullet \ 3 \arrow[u,"\alpha_6"] \arrow[r,swap,"\alpha_3"] & \bullet \ 4 \arrow[r,swap,"\alpha_4"] & \bullet \ 5 \\
\\
\widetilde{E_7} & & & \bullet \ 8 \\
1 \ \bullet \arrow[r,swap,"\alpha_1"] & \bullet \ 2 \arrow[r,swap,"\alpha_2"] & \bullet \ 3 \arrow[r,swap,"\alpha_3"] & \bullet \ 4 \arrow[u,"\alpha_8"] \arrow[r,swap,"\alpha_4"] & \bullet \ 5 \arrow[r,swap,"\alpha_5"] & \bullet \ 6 \arrow[r,swap,"\alpha_6"] & \bullet \ 7
\\
\widetilde{E_8} & & & & & \bullet \ 9 \\
1 \ \bullet \arrow[r,swap,"\alpha_1"] & \bullet \ 2 \arrow[r,swap,"\alpha_2"] & \bullet \ 3 \arrow[r,swap,"\alpha_3"] & \bullet \ 4 \arrow[r,swap,"\alpha_4"] & \bullet \ 5 \arrow[r,swap,"\alpha_5"] & \bullet \ 6 \arrow[u,"\alpha_9"] \arrow[r,swap,"\alpha_6"] & \bullet \ 7 \arrow[r,swap,"\alpha_7"] & \bullet \ 8
\end{tikzcd}
\]
In order to produce the preprojective algebras, one has to consider the doubles, construct the path algebras and factor out one relation. As before, successively multiplying the defining relation by idempotents from either side discloses the following relations.
\begin{align*}
	\widetilde{E_6}: \ \ \ \ \ &\alpha_1 \alpha_1^\ast = \alpha_4^\ast \alpha_4 = \alpha_7^\ast \alpha_7 = 0, \\
	&\alpha_3 \alpha_3^\ast + \alpha_6 \alpha_6^\ast = \alpha_2^\ast \alpha_2, \ \ \ \ \ \ \text{ and } \\
	&\alpha_i \alpha_i^\ast = \alpha_{i-1}^\ast \alpha_{i-1}, \ \ \ \ \ \ \ i = 2,4,7. \\
	\widetilde{E_7}: \ \ \ \ \ &\alpha_1 \alpha_1^\ast = \alpha_6^\ast \alpha_6 = \alpha_8^\ast \alpha_8 = 0, \\
	&\alpha_4 \alpha_4^\ast + \alpha_8 \alpha_8^\ast = \alpha_3^\ast \alpha_3, \ \ \ \ \ \ \text{ and } \\
	&\alpha_i \alpha_i^\ast = \alpha_{i-1}^\ast \alpha_{i-1}, \ \ \ \ \ \ \ i = 2,3,5,6. \\
	\widetilde{E_8}: \ \ \ \ \ &\alpha_1 \alpha_1^\ast = \alpha_7^\ast \alpha_7 = \alpha_9^\ast \alpha_9 = 0, \\
	&\alpha_6 \alpha_6^\ast + \alpha_9 \alpha_9^\ast = \alpha_5^\ast \alpha_5, \ \ \ \ \ \ \text{ and } \\
	&\alpha_i \alpha_i^\ast = \alpha_{i-1}^\ast \alpha_{i-1}, \ \ \ \ \ \ \ i = 2,3,4,5,7.
\end{align*}
As noticed for $\widetilde{A_{n-1}}$ and $\widetilde{D_{n-1}}$, the relations force again $c_1t_1 = c_it_i$ for all $i$ as otherwise $g$ is not an automorphism. Conversely, any nonzero scalars $c_i$, $t_i$ satisfying
$c_it_i = c_1t_1$ give rise to a graded automorphism of $A$.

\begin{proposition}\label{prop51}
Let $Q = \widetilde{E_{m}}$ for $m = 6,7,8$ and $A = \Pi(Q)$. Fix $g$ to be a graded automorphism of $A$ with $g(\alpha_i) = c_i \alpha_i$ and $g(\alpha_i^\ast) = t_i \alpha_i^\ast$ for all $i$. Then $\TrV(g,t)$ equals
\begin{align*}
\begin{pmatrix} 1 + c_1t_1 t^2 & -c_1t &  & & & \\[1ex] -t_1t & 1+c_1t_1t^2 & -c_2t & & & \\[1ex] & -t_2t & 1+c_1t_1t^2 & -c_3t & & -c_6t & \\[1ex] & & -t_3 t & 1+c_1t_1t^2 & -c_4t & & \\[1ex] & & & -t_4t & 1+c_1t_1t^2 & & \\[1ex] & & -t_6t & & & 1+c_1t_1t^2 & -c_7t \\[1ex] & & & & & -t_7t & 1+c_1t_1t^2 \end{pmatrix}^{-1} \cdot \begin{pmatrix} 1 \\[1ex] 1 \\[1ex] 1 \\[1ex] 1 \\[1ex] 1 \\[1ex] 1 \\[1ex] 1 \end{pmatrix}
\end{align*}
for $\widetilde{E_6}$,
\begin{small}
\begin{align*}
\begin{pmatrix} 1 + c_1t_1 t^2 & -c_1t & & & & & & \\[1ex] -t_1t & 1+c_1t_1t^2 & -c_2t & & & & & \\[1ex] & \ddots & \ddots & \ddots & & & & \\[1ex] & & -t_3t & 1+c_1t_1t^2 & -c_4t & & & -c_8t \\[1ex] & & & -t_4t & 1+c_1t_1t^2  & -c_5t & & \\[1ex] & & & & -t_5t & 1+c_1t_1t^2 & -c_6t & \\[1ex] & & & & & -t_6t & 1+c_1t_1t^2 & \\[1ex] & & & -t_8t & & & & 1+c_1t_1t^2 \end{pmatrix}^{-1} \cdot \begin{pmatrix} 1 \\[1ex]  1 \\[1ex] \vdots \\[1ex] 1 \\[1ex] 1 \end{pmatrix}
\end{align*}
\end{small}
for $\widetilde{E_7}$ and
\begin{align*}
\begin{pmatrix} 1 + c_1t_1 t^2 & -c_1t & & & & \\[1ex] -t_1t & 1+c_1t_1t^2 & -c_2t & & &  \\[1ex] & \ddots & \ddots & \ddots & & & \\[1ex] & & -t_5t & 1+c_1t_1t^2 & -c_6t & & -c_9t \\[1ex] & & & -t_6t & 1+c_1t_1t^2  & -c_7t & \\[1ex] & & & & -t_7t & 1+c_1t_1t^2 & \\[1ex]& & & -t_9t & & & 1+c_1t_1t^2  \end{pmatrix}^{-1} \cdot \begin{pmatrix} 1 \\[1ex]  1 \\[1ex] \vdots \\[1ex] 1 \\[1ex] 1 \end{pmatrix}
\end{align*}
for $\widetilde{E_8}$. Similar formulas arise for $\TrVL(g,t)$ where the matrices are replaced by their transposes.

\begin{proof}
This proof is similar to the proof of Proposition \ref{prop42}. Again, we need to find the commutative diagrams whose rows are projective resolutions of the simple modules. They are immediate from knowing the generators, their target vertices and the relations starting at $e_i$ as done before. Then, the calculations for the connecting maps are the same as the corresponding ones in Proposition \ref{prop42}. For instance, for $\widetilde{E_6}$, the idempotent $e_1$ behaves like $e_1$ in $\widetilde{D_{n-1}}$, $e_5$ and $e_7$ behave like $e_n$, the formula for $e_3$ can be extracted from $e_{n-2}$ and at last the idempotents $e_2, e_4$ and $e_6$ fall into the generic case for $\widetilde{D_{n-1}}$. Similar equations to \eqref{eq3new} arise and written in matrix form together with remembering that $\Tr_{S_i}(\tilde{g}_i,t) = 1$ finishes the claim for $\TrV(g,t)$. The case $\TrVL(g,t)$ follows from analogous calculations.
\end{proof}
\end{proposition}

\begin{remark}
It is no coincidence that the matrices appearing in $\TrV(g,t)$ and $\TrVL(g,t)$ are transposes of each other in all cases. The matrix $M(c_i,t_i,t)$ is the inverse of the matrix trace of $g$ but to prove this in details would go beyond the scope of this paper.
\end{remark}

\section{A homological determinant for preprojective algebras}\label{section7}
Let $A = \Pi(Q)$ be an algebra as in Convention \ref{con214} and fix a graded automorphism $g$ of $A$ by $g(\alpha_i) = c_i \alpha_i$ and $g(\alpha_i^\ast) = t_i \alpha_i^\ast$ for some nonzero scalars $c_i$, $t_i \in k$. This section works towards Corollary \ref{cor77} which proves that $\Hdet(g) = c_1t_1 I_n$ ($n = m+1$ for $\widetilde{E_m}$). 

\begin{proposition}\label{prop71}
The following reciprocity statement holds:
\[
	\Tr_{Ae_i}\left(\left(g|_{Ae_i}\right)^{-1},t^{-1}\right) = c_1t_1 t^2 \Tr_{e_iA}(g|_{e_iA},t).
\]
\begin{proof}
We are going to use the formulas for $\TrV(g,t)$ and $\TrVL(g,t)$ from the Propositions \ref{prop33}, \ref{prop42} and \ref{prop51}. There, equations of the form
\[
	\TrV(g,t) = M(c_i,t_i,t)^{-1} \cdot \begin{pmatrix} 1 \\ \vdots \\ 1 \end{pmatrix} \ \ \ \text{ and } \ \ \ \TrVL(g,t) = M(t_i,c_i,t)^{-1} \cdot \begin{pmatrix} 1 \\ \vdots \\ 1 \end{pmatrix} = M(c_i,t_i,t)^{-T} \cdot \begin{pmatrix} 1 \\ \vdots \\ 1 \end{pmatrix} 
\]
with explicit descriptions were found. An easy calculation shows 
\[
	\left((c_1t_1)t^2\right)^{-1} M(c_i,t_i,t) = M(t_i^{-1},c_i^{-1},t^{-1})
\]
which is the matrix of the formula corresponding to $\TrVL(g^{-1},t^{-1})$. It follows that $c_1t_1 t^2 \cdot \TrV(g,t) = \TrVL(g^{-1},t^{-1})$ and looking at the $i$-th entry of the resulting vectors on each side reveals the claim.
\end{proof}
\end{proposition}

The next definition is taken from \cite[Definition 1.3]{MR1758250}.

\begin{definition}
Let $R$ be an $\NN$-graded locally finite $k$-algebra and call $\mathfrak{u} = R_{\geq 1}$. Let $M$ be an $R$-module such that $H^i_{\mathfrak{u}}(M) = 0$ for large $i$ and each $H^i_{\mathfrak{u}}(M)$ is locally finite and right bounded. Let $g \in \Aut_{\gr}(R)$ and $f: M \to M$ be a $g$-linear homomorphism. We define:
\[
	\BTr_M^R(f,t) = \sum_i{(-1)^i\Tr_{H^i_{\mathfrak{u}}(M)}(H^i_{\mathfrak{u}}(f),t)}.
\]
If $M$ is left bounded and locally finite, we call $f$ \textit{rational} over $k$ if the following two conditions are satisfied:
\begin{itemize}
	\item[(1)]	$\Tr_M(f,t)$ and $\BTr_M^R(f,t)$ are rational functions over $k$ (inside $k((t))$ and $k((t^{-1}))$, respectively) and \\[1ex]
	\item[(2)] $\Tr_M(f,t) = \BTr_M^R(f,t)$ as rational functions over $k$.
\end{itemize}
If there is no confusion about the algebra $R$, we may omit the superscript.
\end{definition}

Using the reciprocity property from Proposition \ref{prop71}, we find a connection between $\BTr_{e_iA}(g_i,t)$ and the corresponding trace function: 

\begin{proposition}\label{prop73}
Let $g \in \Aut_{\gr}(A)$. The projective $A$-module $e_iA$ satisfies
\[
	\BTr_{e_iA}^A(g,t) = \Hdet(g)_{ii}^{-1} c_1t_1 \Tr_{e_iA}(g,t)
\]
as rational functions over $k$.
\begin{proof}
From a dual version of \cite[Lemma 2.8]{DanToDo} we know that 
\begin{align*}
	H^2_{\mathfrak{m}}(e_iA) &= \varinjlim{\Ext^2(A/A_{\geq m},e_iA)} = e_i \varinjlim{\Ext^2(A/A_{\geq m},A)} = e_i A^\vee[2].
\end{align*}
Since $H^j_{\mathfrak{m}}(-)$ is an additive functor and $A = \bigoplus_i{e_iA}$, Lemma \ref{lemma63} gives that $H^j_{\mathfrak{m}}(e_iA)$ is nonzero only if $j = 2$. This allows the calculation as in \cite[Proposition 4.2]{MR1758250}:
\begin{align*}
	\BTr_{e_iA}(g,t) &= \Tr_{H^2_{\mathfrak{m}}(e_iA)}(H^2_{\mathfrak{m}}(g),t) \\
	&= \Tr_{e_i A^\vee[2]}(\Hdet(g)_{ii}^{-1} \left(g^{-1}\right)^\ast,t)
	 \\
	&= \Hdet(g)_{ii}^{-1}t^{-2} \Tr_{e_i A^\vee}\left( \left(g^{-1}\right)^\ast,t\right) \\
	&= \Hdet(g)_{ii}^{-1}t^{-2} \Tr_{Ae_i}\left(g^{-1},t^{-1}\right) \\
	&= \Hdet(g)_{ii}^{-1}t^{-2} c_1t_1 t^2 \Tr_{e_iA}\left(g,t\right) \\
	&= \Hdet(g)_{ii}^{-1} c_1t_1 \Tr_{e_iA}\left(g,t\right)
\end{align*}
where $H^2_{\mathfrak{m}}(g|_{e_iA}) = \Hdet(g)_{ii}^{-1}\left(g^{-1}\right)^\ast$ comes from the construction of $\Hdet$ in Lemma \ref{lemma65}. The second to the last equality was proved in Proposition \ref{prop71}.
\end{proof}
\end{proposition}

Before we can proceed let's recall the definition of twisted Calabi-Yau algebras. By \cite[Theorem 5.15]{DanToDo}, these algebras form a not necessarily connected generalization of Artin-Schelter regular algebras.

\begin{definition}[{\cite[Definition 1.2]{DanToDo}}]\label{CYdef}
Let $B$ be a $k$-algebra and call $B^e = B \otimes_k B^{\op}$. We say that $B$ is \textit{twisted Calabi-Yau of dimension $d$} if
\begin{itemize}[topsep=0pt, partopsep=0pt,leftmargin=0.75in]
	\item[(CY 1)] $B$ has a resolution of finite length by finitely generated projective $(B,B)$-bimodules and
	\item[(CY 2)] there exists an invertible $k$-central $(B,B)$-bimodule $U$ such that
	\[
		\Ext_{B^e}^i(B,B^e) \cong \begin{cases} 0, & i \neq d \\ U, & i = d \end{cases}
	\]
	as $(B,B)$-bimodules, where each $\Ext_{B^e}^i(B,B^e)$ is considered as a right $B^e$-module via the ``inner" right $B^e$-structure of $B^e$.
\end{itemize}
\end{definition}

\begin{theorem}\label{thm75}
Let $B$ be a noetherian $\NN$-graded locally finite twisted Calabi-Yau algebra of global dimension $d \geq 1$ with degree zero piece $B_0 \cong k^n$. Denote the primitive orthogonal idempotents by $e_i$. Further assume that there is a noetherian subring $S$ of $B$ which is Artin-Schelter regular of global dimension $d$ such that $_SB$ and $B_S$ are finitely generated. Let $g$ be a graded automorphism of $B$ such that $g(e_i) = e_i$ and $g(S) = S$. Let $f: e_iB \to e_iB$ be a graded $g$-linear automorphism. Then $\BTr_{e_iB}(f,t) = \Tr_{e_iB}(f,t)$. In other words, $f$ is rational over $k$.
\begin{proof}
The algebra $B$ is automatically a generalized Gorenstein algebra by \cite[Theorem 5.2]{DanToDo} and \cite[Theorem 5.15]{DanToDo}. As before, let $\mathfrak{r} = B_{\geq 1}$ and $\mathfrak{t} = S_{\geq 1}$. Then, by \cite[Lemma (2.8)]{OtherPaper} and \cite[Proposition (2.10)]{OtherPaper}, $B$ satisfies the $\chi$-condition. In addition, $S$ satisfies the $\chi$-condition by \cite[Theorem 8.1(1)]{MR1304753}. We recall from the proof of Lemma \ref{lemma67} that \cite[Proposition 7.2(2)]{MR1304753} and \cite[Theorem 8.3(3)]{MR1304753} imply $H_{\mathfrak{r}}^j(M) = H_{\mathfrak{t}}^j(M)$ for all graded right $S$-modules $M$ and all $j \geq 1$. This isomorphism is functorial as described in \cite[Lemma 4.4($2^\circ$)]{MR1758250}, and so $f: e_iB \to e_iB$ induces the same maps on $H_{\mathfrak{r}}^j(e_iB)$ and $H_{\mathfrak{t}}^j(e_iB)$ for $j \geq 1$. By Lemma \ref{lemma63}, $H^j_{\mathfrak{r}}(e_iB)$ and $H^j_{\mathfrak{s}}(e_iB)$ are nonzero if and only if $j = N$. Therefore,
\[
\BTr_{e_iB}^B(f,t) = (-1)^N\Tr_{H^N_{\mathfrak{r}}(e_iB)}(H_{\mathfrak{r}}^N(f),t) = (-1)^N \Tr_{H_{\mathfrak{t}}^N(e_iB)}(H_{\mathfrak{t}}^N(f),t) = \BTr_{e_iB}^S(f,t).
\]
Since for the trace function only the vector space structure matters, $\Tr_{e_iB}(f,t)$ does not depend on whether $e_iB$ is understood as an $S$- or $B$-module. Now, $S$ is Artin-Schelter regular and \cite[Proposition 4.2]{MR1758250} says that $f$ is rational over $k$, i.e. 
\[
	\BTr_{e_iB}^{S}(f,t) = \Tr_{e_iB}(f,t)
\]
as $S$-modules. The claim follows by combining the two equations.
\end{proof}
\end{theorem}

\begin{remark}
For $Q = \widetilde{A_{n-1}}$, $n \geq 2$, it is possible to concretely find a polynomial ring $S$ in two variables over which $A = \Pi(Q)$ is finitely generated which allows one to apply Theorem \ref{thm75} to $\Pi(Q)$. This requires to introduce some notation. In \cite{OtherPaper} it was shown that the following definition is well-defined. The \textit{type} $\type(\beta) = (m,n)$ of a path $\beta \in A$ is defined by $m = \#$ of nonstar arrows in $\beta$ and $n = \#$ of star arrows in $\beta$. Two paths $\beta_1$ and $\beta_2$ are \textit{of the same type} if $\type(\beta_1) = \type(\beta_2)$.\\ \indent
Define $\beta_i = \alpha_i \alpha_{i+1} \cdots \alpha_n \alpha_1 \cdots \alpha_{i-1}$ and $\beta_i^\ast = \alpha_{i-1}^\ast \alpha_{i-2}^\ast \cdots \alpha_1^\ast \alpha_n^\ast \cdots \alpha_i^\ast$. It is easy to see that $\beta_i$ and $\beta_i^\ast$ are loops at vertex $i$. Further, we call
\begin{align*}
	\beta = \sum_{i=1}^{n}{\beta_i}, \ \ \ \ \ \ \ \ \ \ \beta^\ast = \sum_{i=1}^{n}{\beta_i^\ast} \ \ \ \ \ \text{ and } \ \ \ \ \ e = \sum_{i=1}^{n}{e_i} = 1.
\end{align*}
One can check that $\beta$ and $\beta^\ast$ commute. We are going to prove that $A$ is a finitely generated module over the polynomial ring $k[\beta,\beta^\ast]$ where $k = k \cdot e$. Let
\[ 
	\BBB = \{\gamma \text{ path in } \bar{Q} \mid \type(\gamma) \leq (n-1,n-1)\}
\]
with respect to the usual product order. Notice that $\BBB$ is a finite set. We recall that an element of the form $\prod_{j}{\alpha_{i_j}}\prod_{v}{\alpha_{u_v}^\ast}$ is called a simple path. Let $x$ be the unique simple path in $\bar{Q}$ from $e_\ell$ to $e_m$ of type $(a,b)$ where $a = \lambda n + r$ and $b = \lambda^\ast n + r^\ast$ for $0 \leq r,r^\ast < n$ and $\lambda, \lambda^\ast \in \NN_0$. Then we have the following
\[
	\beta^\lambda \left(\beta^\ast\right)^{\lambda^\ast} \tilde{x} = \beta_\ell^\lambda \left(\beta_\ell^\ast\right)^{\lambda^\ast} \tilde{x} = x
\]
where $\tilde{x}$ is the unique simple path from $e_\ell$ to $e_m$ of type $(r,r^\ast)$. This implies that $\tilde{x} \in \BBB$. Therefore, $x$ lies in the module generated by $\BBB$ over $k[\beta,\beta^\ast]$. This shows that $A$ is finitely generated over $k[\beta,\beta^\ast]$ which is an Artin-Schelter regular algebra of global dimension $2$ as used in Theorem \ref{thm75}.
\end{remark}

The previous remark shows that Theorem \ref{thm75} may be applied to a preprojective algebra of type $\widetilde{A_{n-1}}$ for $n \geq 2$. Seeking to apply Theorem \ref{thm75} to other preprojective algebras, we have the following result:

\begin{proposition}\label{propend}
Let $A = \Pi(Q)$ be as in Convention \ref{con214}. For some graded automorphism $g \in \Aut_{\gr}(A)$ which fixes the idempotents, $A$ has a subalgebra $S$ isomorphic to a polynomial ring in two variables such that $_SA$ and $A_S$ are finitely generated and $g(S) = S$.
\begin{proof}
We know that the center of $A$ equals $Z(A) = \mathcal{O}(\mathbb{A}^2)^\Gamma$ by \cite[Theorem 1.9.15(2)]{DanNotes}. As a fixed ring of a polynomial ring in two variables, it is a finitely generated noetherian algebra. Since $g$ is surjective we must have that $g(Z(A)) = Z(A)$. Let $F = Z(A)^{\langle g \rangle}$, the fixed ring of $Z(A)$ under the action of $g$ restricted to $Z(A)$. We recall that $A$ is a finitely generated noetherian $k$-algebra. Moreover, since $A$ is a prime polynomial identity algebra (see Proposition \ref{prop25}), \cite[Corollary B.10.6]{drensky2004polynomial} says that $A$ is a finitely generated module over $Z(A)$. Further, $Z(A)$ is a finitely generated $F$-module by Proposition \ref{prop212}. Together these observations give that $A$ is also a finitely generated $F$-module. Then, the Artin-Tate Lemma \cite[Corollary B.11.4]{drensky2004polynomial} applied to the chain
\[
	k \subseteq F = Z(A)^{\langle g \rangle} \subseteq A
\]
gives that $F$ is a finitely generated $k$-algebra. Now, the $k$-transcendence degree of the quotient field of the center of $A$ equals $2$ by \cite[6.10 Proposition]{MR876985} which passes on to $F$ since $Z(A)$ is a finitely generated $F$-module. The Noether Normalization Lemma gives us a subalgebra $S = k[x_1,x_2]$ of $F$ over which $F$ is a finitely generated module. Consequently, we have the chain of algebras
\[
	S = k[x_1,x_2] \subseteq F \subseteq Z(A) \subseteq A
\]
and we showed that each algebra is a finitely generated module over the previous one. This proves that $A$ is a finitely generated $S$-module. In addition to that, $g(S) = S$ since $S$ is a subalgebra of the fixed ring $F = Z(A)^{\langle g \rangle}$.
\end{proof}
\end{proposition}

\begin{corollary}\label{cor77}
Let $Q$ be the quiver corresponding to an extended Dynkin graph of type $\widetilde{A_{n-1}}$ for $n \geq 2$, $\widetilde{D_{n-1}}$ for $n \geq 5$ or $\widetilde{E_m}$ for $m = 6,7,8$ and $A = \Pi(Q)$. Let $g$ be a graded automorphism of $A$ with $g(\alpha_i) = c_i\alpha_i$ and $g(\alpha_i^\ast) = t_i\alpha_i^\ast$. Then $\Hdet(g) = c_1t_1 I_n$ ($n = m+1$ for $\widetilde{E_m}$).
\begin{proof}
Proposition \ref{propend} gives us a subalgebra $S = k[x_1,x_2]$ of $A$ such that $_SA$ and $A_S$ are finitely generated and $g(S) = S$. Clearly, $S$ is Artin-Schelter regular of global dimension $2$. We note that $A$ is a noetherian $\NN$-graded twisted Calabi-Yau algebra of global dimension $2$ by \cite[Proposition (2.11)]{OtherPaper}. Therefore $f = g_i: e_iB \to e_iB$ is rational over $k$ by Theorem \ref{thm75}. Combined with Proposition \ref{prop73} we get
\[
	\Tr_{e_iA}(g_i,t) = \BTr_{e_iA}^A(g_i,t) = \Hdet(g)^{-1}_{ii} c_1t_1 \Tr_{e_iA}(g,t).
\]
Since $g_i(e_i) = e_i$ it is guaranteed that $\Tr_{e_iA}(g_i,t)$ is nonzero and thus $\Hdet(g)_{ii} = c_1t_1$.
\end{proof}
\end{corollary}

\begin{remark}
This paper generalizes the work of Peter J{\o}rgensen and James J. Zhang in \cite{MR1758250} to the non-connected setting while imposing the restriction that the group action fixes the primitive orthogonal idempotents. It would be very interesting to study graded automorphisms which permute the primitive orthogonal idempotents. The author strongly believes that the results are still true in a similar way.  
\end{remark}

Concluding this paper we take a look at an example to see how the theorems apply. 

\begin{example}
Let $A = \Pi(Q)$ for $Q = \widetilde{A_2}$. We define $g$ to map:
\begin{center}
\begin{tabular}{lllllll}
	$\alpha_1 \mapsto \alpha_1$, & & & $\alpha_2 \mapsto -\alpha_2$, & & & $\alpha_3 \mapsto -\alpha_3$, \\
	$\alpha_1^\ast \mapsto \alpha_1^\ast$, & & & $\alpha_2^\ast \mapsto -\alpha_2^\ast$, & & & $\alpha_3^\ast \mapsto -\alpha_3^\ast$
\end{tabular}
\end{center}
and extend the action linearly. Corollary \ref{cor77} guarantees that $\Hdet(g) = I_3$. Clearly, $g$ has order $2$. Using Proposition \ref{prop33} to calculate $\TrV(g,t)$ and summing up the entries, we find the trace of $g$. Further, Molien's Theorem \ref{lemma211} provides the Hilbert series of $A^G$ as follows:
\begin{align*}
\Tr(\id,t) = \frac{3}{(1-t)^2}, \hspace{0.7cm}
\Tr(g,t) = \frac{3-5t+3t^2}{(1-t)(1-t^3)},  \hspace{0.7cm} \text{and} \hspace{0.7cm} \HS_{A^G}(t) = \frac{3-t+3t^2}{(1-t)(1-t^3)}.
\end{align*}
In order to find all generators of the fixed ring $A^G$ for $G = \langle g \rangle$ it is enough to look at simple representations of paths of $A$. A straightforward analysis reveals the following generators
\[
	\alpha_1, \ \ \alpha_1^\ast, \ \ x = \alpha_2 \alpha_3, \ \ y = \alpha_3^\ast \alpha_2^\ast, \ \ u = \alpha_3 \alpha_1 \alpha_2, \ \ v = \alpha_3^\ast \alpha_2^\ast \alpha_1^\ast, \ \ w = \alpha_3 \alpha_3^\ast
\]
of $A^G$. This means our initial guess for a presentation of $A^G$ is the path algebra $B$ of
	
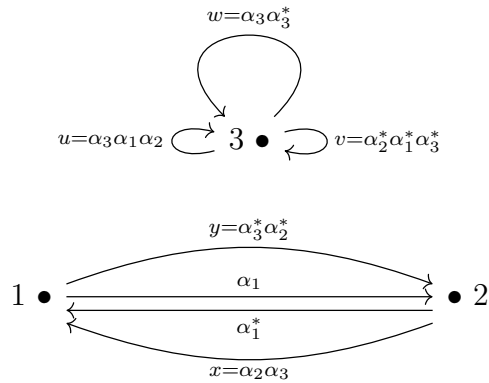
\[
\begin{tikzcd}
& & 3 \ \bullet \arrow[loop right, "v = \alpha_2^\ast \alpha_1^\ast \alpha_3^\ast"] \arrow[loop left, "u =\alpha_3 \alpha_1 \alpha_2"] \arrow[swap,loop, "w = \alpha_3 \alpha_3^\ast"]  \\
& \\
1 \ \bullet \arrow[rrrr,"\alpha_1"] \arrow[rrrr,"y = \alpha_3^\ast \alpha_2^\ast", out = 20, in = 160] & & & & \bullet \ 2 \arrow[llll,"\alpha_1^\ast", shift left = 1ex] \arrow[llll,"x = \alpha_2\alpha_3", out = 200, in = 340, shift left = 1ex] 
\end{tikzcd}
	\]
subject to the following relations suggested by the relations of $A$:
\begin{align*}
yx = (\alpha_1 \alpha_1^\ast)^2, \ \ \  y\alpha_1^\ast \alpha_1 = \alpha_1 \alpha_1^\ast y, \ & \ \ xy = (\alpha_1^\ast \alpha_1)^2, \ \ \ x \alpha_1 \alpha_1^\ast = \alpha_1^\ast \alpha_1 x, \\
w^3 = uv, \ \ \ vu = uv, \ & \ \ wu = uw, \ \ \ wv = vw.
\end{align*}
It is easy to see that the $e_3Be_3$ component is isomorphic to $k[u,v,w]/(w^3-uv)$. In regards to the other connected component one needs to check that the overlaps $yxy, \ yx\alpha_1 \alpha_1^\ast, \ xyx$ and $xy\alpha_1^\ast \alpha_1$ resolve. Then, a Gr\"obner basis calculation, counting words and comparing to $\HS_{A^G}(t)$ gives that this is indeed a presentation for $A^G$.\\ \indent
With Theorem \ref{thm69} in mind, we know that $A^G$ is generalized Gorenstein. To show that this is not obvious, we verify the generalized Gorenstein condition manually. It is well-known that $k[u,v,w]/(w^3-uv)$ is Artin-Schelter Gorenstein. Hence, we only need to look at the other connected component.\\ \indent
By symmetry, it is enough to calculate the minimal projective resolution $P_\bullet \to S_1$ of the simple module $S_1 = e_1A^G/(e_1A^G)_{\geq 1}$. As always, denote the indecomposable projective modules of $A^G$ by $P_i = e_iA^G$ for $i = 1,2$. Then, a tedious calculation verifies that
	\[
\begin{tikzcd}[row sep = 30pt, column sep = 30pt]
\ldots \arrow[r] & P_1 \oplus P_2 \arrow[rrrr, "\begin{pmatrix} \alpha_1 \alpha_1^\ast \ \ \ \ \ \ y \\ \ -x \ \ \ \  -\alpha_1^\ast \alpha_1 \end{pmatrix}  \boldsymbol{\cdot}"] & & & & P_1 \oplus P_2 \arrow[rrr, "\begin{pmatrix} \alpha_1 \alpha_1^\ast \ \ \ \ \ \ y \\ \ -x \ \ \ \  -\alpha_1^\ast \alpha_1 \end{pmatrix}  \boldsymbol{\cdot}"] & & & \ldots
\\
& P_1 \oplus P_2 \arrow[rrrr, "\begin{pmatrix} \alpha_1^\ast \alpha_1 \alpha_1^\ast \ \ \ \alpha_1^\ast y \\ \ -x \ \ \ \ \  -\alpha_1^\ast \alpha_1 \end{pmatrix} \boldsymbol{\cdot}"] & & & & P_2 \oplus P_2 \arrow[rr, "(\alpha_1 \ \ y)  \boldsymbol{\cdot}"] & & P_1\arrow[r] & S_1 \arrow[r] & 0
\end{tikzcd}
\]
is a minimal projective resolution of $S_1$. Applying $\GrHom_{A^G}(-,A^G)$ to $P_{\bullet} \to 0$ and using {\cite[Lemma (2.5)]{OtherPaper}} returns a sequence exact everywhere but at the $2$nd position. A calculation shows that $\GrExt_{A^G}^2(S_1,A^G) \cong S_2^\vee = A^Ge_2/(A^Ge_2)_{\geq 1}$ and by symmetry we also get $\GrExt_{A^G}^2(S_2,A^G) \cong S_1^\vee$. 
\end{example}

\bibliographystyle{amsalpha}
\bibliography{Generalized_Gorensteinness_and_a_homological_determinant_for_preprojective_algebras}

\newcommand{\etalchar}[1]{$^{#1}$}
\def\cprime{$'$}
\providecommand{\bysame}{\leavevmode\hbox to3em{\hrulefill}\thinspace}
\providecommand{\MR}{\relax\ifhmode\unskip\space\fi MR }
\providecommand{\MRhref}[2]{%
  \href{http://www.ams.org/mathscinet-getitem?mr=#1}{#2}
}
\providecommand{\href}[2]{#2}
\begin{thebibliography}{CKWZ16}

\bibitem[ASS06]{MR2197389}
I.~Assem, D.~Simson, and A.~Skowro{\'n}ski, \emph{Elements of the
  representation theory of associative algebras. {V}ol. 1}, London Mathematical
  Society Student Texts, vol.~65, Cambridge University Press, Cambridge, 2006,
  Techniques of representation theory. \MR{2197389}

\bibitem[AZ94]{MR1304753}
M.~Artin and J.~J. Zhang, \emph{Noncommutative projective schemes}, Adv. Math.
  \textbf{109} (1994), no.~2, 228--287. \MR{1304753}

\bibitem[BBK02]{MR1930968}
S.~Brenner, M.~C.~R. Butler, and A.~D. King, \emph{Periodic algebras which are
  almost {K}oszul}, Algebr. Represent. Theory \textbf{5} (2002), no.~4,
  331--367. \MR{1930968}

\bibitem[BGL87]{MR876985}
D.~Baer, W.~Geigle, and H.~Lenzing, \emph{The preprojective algebra of a tame
  hereditary {A}rtin algebra}, Comm. Algebra \textbf{15} (1987), no.~1-2,
  425--457. \MR{876985}

\bibitem[BRS{\etalchar{+}}16]{DanNotes}
G.~Bellamy, D.~Rogalski, T.~Schedler, J.~T. Stafford, and M.~Wemyss,
  \emph{Noncommutative algebraic geometry}, Mathematical Sciences Research
  Institute Publications, Cambridge University Press, 2016.

\bibitem[CKWZ16]{MR3552496}
K.~Chan, E.~Kirkman, C.~Walton, and J.~J. Zhang, \emph{Quantum binary
  polyhedral groups and their actions on quantum planes}, J. Reine Angew. Math.
  \textbf{719} (2016), 211--252. \MR{3552496}

\bibitem[DF12]{drensky2004polynomial}
V.~Drensky and E.~Formanek, \emph{Polynomial identity rings}, Birkh\"auser,
  2012.

\bibitem[DR81]{MR624903}
V.~Dlab and C.~M. Ringel, \emph{Eigenvalues of {C}oxeter transformations and
  the {G}el\cprime fand-{K}irillov dimension of the preprojective algebras},
  Proc. Amer. Math. Soc. \textbf{83} (1981), no.~2, 228--232. \MR{624903}

\bibitem[EE07]{MR2335985}
P.~Etingof and C.-H. Eu, \emph{Koszulity and the {H}ilbert series of
  preprojective algebras}, Math. Res. Lett. \textbf{14} (2007), no.~4,
  589--596. \MR{2335985 (2008j:16048)}

\bibitem[J{\o}r99]{MR1643863}
P.~J{\o}rgensen, \emph{Non-commutative {C}astelnuovo-{M}umford regularity},
  Math. Proc. Cambridge Philos. Soc. \textbf{125} (1999), no.~2, 203--221.
  \MR{1643863}

\bibitem[JZ97]{MR1438180}
N.~Jing and J.~J. Zhang, \emph{On the trace of graded automorphisms}, J.
  Algebra \textbf{189} (1997), no.~2, 353--376. \MR{1438180 (98f:16029)}

\bibitem[JZ00]{MR1758250}
P.~J{\o}rgensen and J.~J. Zhang, \emph{Gourmet's guide to {G}orensteinness},
  Adv. Math. \textbf{151} (2000), no.~2, 313--345. \MR{1758250}

\bibitem[KKZ08]{MR2434290}
E.~Kirkman, J.~Kuzmanovich, and J.~J. Zha{}ng, \emph{Rigidity of graded regular
  algebras}, Trans. Amer. Math. Soc. \textbf{360} (2008), no.~12, 6331--6369.
  \MR{2434290 (2009e:16065)}

\bibitem[KKZ09]{MR2568355}
E.~Kirkman, J.~Kuzmanovich, and J.~J. Zhan{}g, \emph{Gorenstein subrings of
  invariants under {H}opf algebra actions}, J. Algebra \textbf{322} (2009),
  no.~10, 3640--3669. \MR{2568355}

\bibitem[KKZ10]{MR2601538}
E.~Kirkman, J.~Kuzmanovich, and J.~J. Zhang, \emph{Shephard-{T}odd-{C}hevalley
  theorem for skew polynomial rings}, Algebr. Represent. Theory \textbf{13}
  (2010), no.~2, 127--158. \MR{2601538}

\bibitem[MM11]{MR2770441}
H.~Minamoto and I.~Mori, \emph{The structure of {AS}-{G}orenstein algebras},
  Adv. Math. \textbf{226} (2011), no.~5, 4061--4095. \MR{2770441}

\bibitem[Mon80]{MR590245}
S.~Montgomery, \emph{Fixed rings of finite automorphism groups of associative
  rings}, Lecture Notes in Mathematics, vol. 818, Springer, Berlin, 1980.
  \MR{590245}

\bibitem[RR18a]{DanToDo2}
Manuel~L. Reyes and Daniel Rogalski, \emph{Growth of graded twisted calabi-yau
  algebras}, 2018.

\bibitem[RR18b]{DanToDo}
Manuel~L. Reyes{} and Daniel Rogalski, \emph{A twisted calabi-yau toolkit},
  2018.

\bibitem[RRZ14]{MR3250287}
M.~Reyes, D.~Rogalski, and J.~J. Zhang, \emph{Skew {C}alabi-{Y}au algebras and
  homological identities}, Adv. Math. \textbf{264} (2014), 308--354.
  \MR{3250287}

\bibitem[vdB97]{MR1469646}
M.~van~den Bergh, \emph{Existence theorems for dualizing complexes over
  non-commutative graded and filtered rings}, J. Algebra \textbf{195} (1997),
  no.~2, 662--679. \MR{1469646}

\bibitem[Wei18]{dissertation}
S.~Weispfennin{}g, \emph{Invariant theory of preprojective algebras}, Ph.D.
  thesis, UC San Diego, 2018, ProQuest ID: Weispfenning\_ucsd\_0033D\_17915.
  Merritt ID: ark:/13030/m56x488z.

\bibitem[Wei19]{OtherPaper}
S.~Weispfenning, \emph{Properties of the fixed ring of a preprojective
  algebra}, J. Algebra \textbf{517} (2019), 276--319.

\bibitem[YZ99]{MR1674648}
A.~Yekutieli and J.~J. Zhang, \emph{Rings with {A}uslander dualizing
  complexes}, J. Algebra \textbf{213} (1999), no.~1, 1--51. \MR{1674648}

\end{thebibliography}

\end{document}